\documentclass[smallextended]{svjour3}

\usepackage{algorithmic} 
\usepackage{algorithm}

\usepackage{xcolor}
\usepackage{float} 

\usepackage{amsmath}
\usepackage{amssymb}

\usepackage{hyperref}

\usepackage{physics}

\newcommand{\hess}{\nabla^2}
\newcommand{\lap}{\Delta}

\newcommand{\epsp}{\epsilon^+}
\newcommand{\epsm}{\epsilon^-}
\newcommand{\omgp}{\Omega^+}
\newcommand{\omgm}{\Omega^-}
\newcommand{\R}{\mathbb{R}}
\newcommand{\jump}[1]{\left[#1\right]}

\renewcommand{\norm}[1]{\left\lVert#1\right\rVert}

\newcommand{\hxk}{{\vb{\hat{x}}_k}}
\newcommand{\bxi}{{\vb{x_i}}} 
\newcommand{\hx}{{\vb{\hat{x}}}}
\newcommand{\vn}{{\vb{n}}}

\newcommand{\ujbr}{u_{\vb j \in B_r}}

\RenewDocumentCommand\partialderivative{ s o m g g d() }
{ 
	\IfBooleanTF{#1}
	{\let\fractype\flatfrac}
	{\let\fractype\frac}
	\IfNoValueTF{#4}
	{
		\IfNoValueTF{#6}
		{\fractype{\partial \IfNoValueTF{#2}{}{^{#2}}}{\partial #3\IfNoValueTF{#2}{}{^{#2}}}}
		{\fractype{\partial \IfNoValueTF{#2}{}{^{#2}}}{\partial #3\IfNoValueTF{#2}{}{^{#2}}} \argopen(#6\argclose)}
	}
	{
		\IfNoValueTF{#5}
		{\fractype{\partial \IfNoValueTF{#2}{}{^{#2}} #3}{\partial #4\IfNoValueTF{#2}{}{^{#2}}}}
		{
			\IfNoValueTF{#2}
			{\fractype{\partial^2 #3}{\partial #4 \partial #5}}
			{\fractype{\partial^{#2} #3}{\partial #4 \partial #5}}
		}
	}
}

\usepackage{mathtools}

\DeclarePairedDelimiter\floor{\lfloor}{\rfloor}

\usepackage{bm}
\usepackage{caption}
\usepackage{subcaption}

\begin{document}

\title{A Compact Coupling Interface Method with Second-Order Gradient Approximation for Elliptic Interface Problems}
\titlerunning{Compact Coupling Interface Method for Elliptic Interface Problems}
\author{Ray Zirui Zhang \and Li-Tien Cheng}

\institute{Ray Zirui Zhang \\
\email{zzirui@ucsd.edu} \\
Li-Tien Cheng \\
\email{lcheng@math.ucsd.edu} \\
\at Department of Mathematics, University of California, San Diego 
\at 9500 Gilman Drive, La Jolla, CA, 92093-0112, USA}


\maketitle


\begin{abstract}

We propose the Compact Coupling Interface Method (CCIM), a finite 
difference method capable of obtaining second-order accurate approximations of not only 
solution values but their gradients, for elliptic complex interface problems with interfacial 
jump conditions. 
Such elliptic interface boundary value problems with interfacial jump conditions are a critical part of numerous applications in fields such as heat conduction, 
fluid flow, materials science, and protein docking, to name a few.  A typical 
example involves the construction of biomolecular shapes, where such elliptic interface
problems are in the form of linearized Poisson-Boltzmann equations, involving 
discontinuous dielectric constants across the interface, that govern electrostatic 
contributions.  Additionally, when interface dynamics are involved, the normal velocity 
of the interface might be comprised of the normal derivatives of solution, which can be approximated to second-order by our method, resulting in accurate interface dynamics.
Our method, which can be formulated in arbitrary spatial dimensions,
combines elements of the highly-regarded Coupling Interface Method, for such elliptic
interface problems, and Smereka's second-order accurate discrete delta function.  The 
result is a variation and hybrid with a more compact stencil than that found in the 
Coupling Interface Method, and with advantages, borne out in numerical experiments 
involving both geometric model problems and complex biomolecular surfaces, in more 
robust error profiles.
\end{abstract}


\section{Introduction}

%
%
%
%
%
%
%

\subsection{Applications}

Elliptic interface problems with interfacial jump conditions can be found at the 
heart of a variety of physical and biological problems involving interfaces.  These 
interfaces may be material interfaces or phase boundaries, static or dynamic, and in 
subjects relating to heat conduction \cite{gibouFourthOrderAccurate2005,gibouSecondOrderAccurateSymmetricDiscretization2002}, fluid dynamics \cite{bochkovSolvingEllipticInterface2020}, materials science \cite{houHybridMethodMoving1997,kafafyThreedimensionalImmersedFinite2005}, electromagnetics \cite{zhaoHighOrderMatched2010,hadleyHighaccuracyFinitedifferenceEquations2002}, or 
electrostatics \cite{shuAccurateGradientApproximation2014,zhouVariationalImplicitSolvation2014,zhongImplicitBoundaryIntegral2018},
tumor growth \cite{macklinEvolvingInterfacesGradients2005,macklinNewGhostCell2008}.
The interfacial jump conditions are due, frequently, to material 
properties and sources that are discontinuous, or have discontinuous derivatives, 
across the interface.  This leads to solutions that also have discontinuities in 
values or derivatives at the interface.  These discontinuities, especially when 
they are large, are what presents the main difficulties in this problem.


Our motivating application that fits into this framework involves biomolecular shapes.  
Consider a set of solute atoms making up a biomolecule, or several biomolecules, such 
as proteins involved in a docking process.  One interest in this situation is how 
the atoms affect the solvent that it resides in, usually a solution resembling 
salt-water.  The implicit solvation approach introduces an interface to separate a
continuously modeled solvent from the solute atoms and vacuum \cite{dzubiellaCouplingHydrophobicityDispersion2006,dzubiellaCouplingNonpolarPolar2006}
It additionally pairs with this a free energy involving contributions such as nonpolar van der Waals 
forces, surface effects, and electrostatics, with the minimizer serving as the desired
interface \cite{wangLevelSetVariationalImplicitSolvent2012,zhouVariationalImplicitSolvation2014,zhangCouplingMonteCarlo2021,zhangBinaryLevelSet2023}.
The electrostatics portion here provides the elliptic interface problem with interfacial jump conditions we are interested in, arising from linearization
approximations of the governing Poisson-Boltzmann equation 
\cite{zhouVariationalImplicitSolvation2014,izzoCorrectedTrapezoidalRuleIBIM2022,zhongImplicitBoundaryIntegral2018,holstMultigridSolutionPoisson1993,holstTreatmentElectrostaticEffects1994}.

\subsection{Setup}

\begin{figure}[!htbp]
\centering
\includegraphics[width=0.5\textwidth]{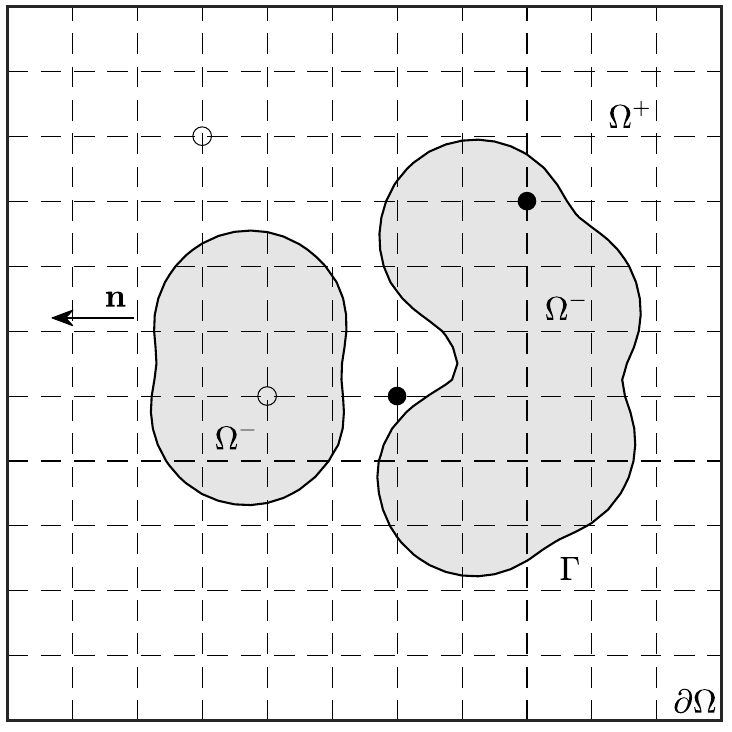}
\caption{Schematic for the elliptic interface problem. $\Gamma$ is an interface that separates a cubical domain $\Omega$ with boundary $\partial \Omega$ into $\omgp$ and $\omgm$. The normal to the interface is denoted as $\vb n$. The dashed lines are the grid lines of the uniform mesh.  The circles are interiors points, where standard central difference can be used.
The disks are examples of on-front points, where special stencils are required. }
\label{f:schematic}
\end{figure}

Let $\Omega\subset\R^d$ be a rectangular box and consider an orientable $\mathcal{C}^1$ hypersurface $\Gamma$ that 
separates it into an inside region $\omgm$ and an outside region $\omgp$.  Also let 
$\vn$ denote outward unit normal vectors on the interface (see 
Fig.~\ref{f:schematic}).  Then, for given functions $\epsilon$, $f$, $a$: 
$\Omega \to \R$, possibly discontinuous across the interface, and given functions 
$\tau$, $\sigma$: $\Gamma \to \R$, our specific elliptic interface problem of 
interest, with interfacial jump conditions, takes the form:
\begin{equation}
\label{eq:pde}
\begin{cases}
-\grad \cdot(\epsilon \grad  u) + a u = f & \text{in }\Omega \setminus \Gamma, \\
\left[u\right] = \tau,\quad \left[\epsilon\grad  u\cdot \vn\right] = \sigma & \text{on }\Gamma, \\
u = g & \text{on }\partial \Omega. \\
\end{cases}
\end{equation}
Here, for any $v:\Omega \to \R$ a function and $x\in\Gamma$, we employ the commonly used 
notation of $[v]$ to denote the jump of $v$ across the interface at $x$:
\begin{equation}
 [v] = v^+ - v^-.
\end{equation}
The superscript $+$ or $-$ denotes the limiting value of a function from $\Omega^+$ or $\Omega^-$, respectively.
Additionally, we will refer to $\epsilon \grad u \cdot \vn$ as the flux;
$\tau$, $\sigma$: $\Gamma \to \R$ as the value of the jump conditions; $\Omega$ as the 
computational domain; and $g$ as the value of the Dirichlet boundary conditions on 
$\partial\Omega$.  Note, we 
allow the dimension $d$ to be general for the formulation of our method, but restrict
to the case $d = 3$ for computations, which we find to have sufficient complexity for 
many real-world problems.

In our motivating application, this form is achieved in a linearized Poisson-Boltzmann
equation for electric potential $u$, charge density $f$, and dielectric coefficient 
$\epsilon$, which can take on values of around $1$ or $2$ in the solute region and $80$ 
in the solvent region \cite{shuAccurateGradientApproximation2014,zhongImplicitBoundaryIntegral2018,zhouVariationalImplicitSolvation2014}.

\subsection{Interface Dynamics and Gradients}

In applications that concern both such elliptic interface problems and moving interfaces, 
the velocity of the interface may depend on the gradient of the solution.  This makes
it imperative to accurately calculate not only the solution's values, but its gradient
as well \cite{macklinEvolvingInterfacesGradients2005,liMinimizationElectrostaticFree2009,zhouVariationalImplicitSolvation2014}. 
In fact, it is preferrable for errors to be measured under the infinity norm, to obtain accurate pointwise velocities that can be used in both front tracking methods and level-set methods for interface dynamics.

In our motivating biomolecular application, gradient descent on the free-energy 
functional can be used to capture the solute-solvent interface of interest 
\cite{dzubiellaCouplingNonpolarPolar2006}. This introduces a time 
variable and interface dynamics into the originally static problem.  Additionally,
gradient descent in this case translates to a normal velocity that depends on the 
effective dielectric boundary force, which in turn depends on the jump of the
normal derivative of the electrostatic potential at the interface
\cite{liMinimizationElectrostaticFree2009,zhouVariationalImplicitSolvation2014}.

In another and perhaps better known example, the Stefan problem 
\cite{gibouSecondOrderAccurateSymmetricDiscretization2002,gibouFourthOrderAccurate2005} is used to model an interface 
separating ice from water, where the interface moves when ice melts or water freezes.  
In simple terms, for a given configuration of the ice and water, an elliptic interface 
problem, with $\epsilon$ as thermal conductivity, can be solved for the temperature $u$.  
The interface between the ice and water then evolves due to this temperature, with its 
normal velocity depending on the jump in the normal derivative of the temperature.

From these considerations, we can finally formulate our goal: to accurately and 
efficiently solving elliptic interface problems with jump conditions for both solution 
and gradient.

\subsection{Finite Difference Methods}\label{ss:fdm}

There are many existing approaches that solve elliptic interface problems with jump 
conditions, however, they can roughly be categorized by the framework they employ.  
An example is the Coupling Interface Method (CIM) 
\cite{chernCouplingInterfaceMethod2007}, which is a finite difference method,
employing a fixed grid that provides ease in resolution and in the construction of 
accurate schemes such as central differencing.  In fact, the main classes of methods 
found include finite difference, finite element, boundary integral, finite volume, and 
deep learning methods, each with their own inherent advantages and disadvantages and
suitability for chosen applications.

An example of such suitability involves elliptic interface problems with their 
jump boundary conditions set on moving boundaries.  In this setting, one has the added 
concern of how well a chosen framework interacts with a chosen numerical representation 
of the boundary interface.  If the level-set method \cite{osherLevelSetMethods2003a,osherFrontsPropagatingCurvaturedependent1988} -- a finite difference method that 
has proven itself in its handling of complicated dynamics, especially topological 
changes -- were used to represent the interface, then a finite difference method for 
the elliptic interface problem would be the natural fit; in fact, both methods could 
share the same underlying grid.  

This example, though, is exactly the reason that we decide to concentrate on finite 
difference methods in this paper.  In our biomolecule application of interest, under 
implicit solvation, where the solvent is represented as a continuous medium, the goal 
is to capture the interface that separates a biomolecule from the surrounding solvent.  
The Variational Implicit Solvation Model (VISM) constructs the free energy of the system 
and seeks a minimizing interface \cite{dzubiellaCouplingNonpolarPolar2006,dzubiellaCouplingHydrophobicityDispersion2006,chengInterfacesHydrophobicInteractions2009}.
One approach to obtaining the minimizer is gradient descent on the energy, resulting in the flow of an initial interface to, in steady state, a minimizer. 
This interface flow was handled, in \cite{chengInterfacesHydrophobicInteractions2009,wangLevelSetVariationalImplicitSolvent2012,zhouVariationalImplicitSolvation2014}, exactly by the level-set 
method.  Thus, to build on the work and results there, in adding electrostatic
effects, which can be described by elliptic interface problems with jump boundary 
conditions on the moving interfaces, we choose the finite difference method as our 
framework.  

Again, this is not to say other frameworks are not, in other circumstances, 
advantageous.  We simply believe a finite difference method would best fit
our problem of interest in our application of interest and likely others as well.  
Under other frameworks, we note especially the boundary integral methods of
\cite{bealeGridbasedBoundaryIntegral2004,zhongImplicitBoundaryIntegral2018,guoSolvingThreedimensionalInterface2021,izzoCorrectedTrapezoidalRuleIBIM2022};
finite volume methods of \cite{bochkovSolvingEllipticInterface2020,guittetSolvingEllipticProblems2015}; finite element methods of
\cite{chenFiniteElementMethods1998,huangMortarElementMethod2002,liNewFormulationsInterface2003,guoGradientRecoveryElliptic2018}, particularly those involving unfitted meshes that are suitable for moving interfaces
\cite{guoFixedMeshMethod2019,guoSolvingParabolicMoving2021,gongImmersedInterfaceFiniteElementMethods2008,beckerNitscheExtendedFinite2009,burmanGhostPenalty2010,chuNewMultiscaleFinite2010,zuninoUnfittedInterfacePenalty2011,barrauRobustVariantNXFEM2012,adjeridEnrichedImmersedFinite2023,massingStabilizedNitscheFictitious2014,guzmanFiniteElementMethod2017a};
and, recently, the deep learning methods of \cite{huDiscontinuityCapturingShallow2021,guoDeepUnfittedNitsche2022}.
This list, furthermore, is by no means an exhaustive account of methods under other 
frameworks.  The topic of elliptic interface 
problems with jump boundary conditions is a very active area of research with a large 
body of literature on the algorithms, analysis, and applications.

As discussed in \cite{chernCouplingInterfaceMethod2007}, within finite difference methods, there are 
essentially three types of approaches: regularization, dimension unsplitting, and 
dimension splitting approaches.  A regularizatation approach applies smoothing 
techniques to discontinuous coefficients, or regularization techniques to singular 
sources \cite{tornbergNumericalApproximationsSingular2004,tornbergRegularizationTechniquesNumerical2003}.  An major example 
of this that is related to our problem of interest is the Immersed Boundary Method (IBM) 
\cite{peskinImmersedBoundaryMethod2002,peskinNUMERICALANALYSISBLOODFLOWHEART1977}.  In dimension unsplitting 
approaches, finite difference methods are derived from local Taylor expansions in 
multi-dimensions.  One popular method in this category is the Immersed Interface Method 
\cite{levequeImmersedInterfaceMethod1994} and its various extensions, including the 
Maximum Principle Preserving Immersed Interface Method (MIIM)\cite{liMaximumPrinciplePreserving2001}, the 
Fast Immersed Interface Method (FIIM) \cite{liFastIterativeAlgorithm1998}, and the Augmented Immersed 
Interface Method (AIIM) \cite{liAccurateSolutionGradient2017}.  For dimension splitting approaches, 
the finite difference methods are derived from Taylor expansions in each dimension. 
This category includes the Ghost Fluid Method 
\cite{liuConvergenceGhostFluid2003,liuBoundaryConditionCapturing2000,fedkiwNonoscillatoryEulerianApproach1999}, the 
Explicit-jump Immersed Interface Method (EJIIM) \cite{wiegmannExplicitjumpImmersedInterface2000}, the 
Decomposed Immersed Interface Method (DIIM) \cite{berthelsenDecomposedImmersedInterface2004}, the 
Matched Interface and Boundary Method (MIB) \cite{zhouHighOrderMatched2006,yuMatchedInterfaceBoundary2007,yuThreedimensionalMatchedInterface2007}, and CIM and its variation, Improved Coupling Interface Method (ICIM) 
\cite{chernCouplingInterfaceMethod2007,shuAccurateGradientApproximation2014,shuAugmentedCouplingInterface2010}.  For a more detailed discussion on these different types, we refer readers to \cite{chernCouplingInterfaceMethod2007}.

In fact, towards our goal, we will be combining two pieces of work set in
the finite difference framework, namely CIM 
\cite{chernCouplingInterfaceMethod2007} and Smereka's work on discrete delta functions 
\cite{smerekaNumericalApproximationDelta2006}.  The former is able to produce 
second-order accurate solutions and first order accurate derivatives of general elliptic 
interface problems with jump boundary conditions, while the latter can perform
second-order accurate operations with a discrete delta function derived within a 
specific class of elliptic interface problems with jump boundary conditions.
To see the connection between delta functions and elliptic interface problems,
consider the Green's function for Laplace's equation: $\lap u = \delta_\Gamma$ with $u = 0$ on $\partial \Omega$, where $\delta_\Gamma$ is the delta function supported on the interface $\Gamma$. This is equivalent to the elliptic interface problem: $\lap u = 0$ with $[u] = 0$ and $[\grad u \cdot \vb n] = 1$.
This connection is used to construct a second-order accurate discrete approximation of the delta function in \cite{smerekaNumericalApproximationDelta2006}.



\subsection{The Coupling Interface Method}\label{ss:cim}

Among these finite difference methods, CIM is one of the top ones in terms of accuracy, 
in both solutions values and gradients; ease of use; and detail of study (see 
\cite{chernCouplingInterfaceMethod2007,shuAccurateGradientApproximation2014,shuAugmentedCouplingInterface2010}).  
In CIM, the standard central differencing stencil is used when the grid point is away from the interface.
When the grid point is next to the interface that forbids the use of the standard central differencing stencil, CIM uses polynomial approximations on either side of the interface, in each dimension, and connects them with jump conditions at the interface.
This leads to a coupled linear system of equations to be solved for the principal second-order derivatives in terms of values at gridpoints.

One version of this approach, called CIM1, chooses linear polynomials and lower-order 
approximations of mixed derivatives for a lower-order but widely applicable approach; 
another, called CIM2, chooses quadratic polynomials and higher-order approximations of 
mixed derivatives for a higher-order approach that, however, requires certain larger 
stencils.  CIM is a hybrid of these that uses CIM2 approximations at points where 
the stencils allow, and CIM1 approximations at all other gridpoints, called exceptional 
points.  Note, these exceptional points do commonly exist but, as noted in 
\cite{chernCouplingInterfaceMethod2007}, not in great numbers, allowing CIM to be second-order 
accurate in solution values under the infinity norm.  For gradients, however, this
approach is only first-order accurate, especially at exceptional points.

The Improved Coupling Interface Method (ICIM) \cite{shuAccurateGradientApproximation2014} fixes this 
issue and achieves uniformly second-order accurate gradient by incorporating two recipes 
that handle exceptional points.  One attempts to ``shift'': at gridpoints where it is 
difficult to achieve a valid first-order approximation of the principal or mixed 
second-order derivative, the finite difference approximations at adjacent gridpoints,
on the same side, are instead shifted over.  The other attempts to ``flip'': at some 
gridpoints, the signature of the domain (inside or outside) can be flipped, so that the 
usual second-order CIM2 discretization may apply, allowing for second-order accurate 
solutions and gradients for the flipped interface.  Extrapolation from neighboring 
nonflipped gridpoints originally on the same side can then be used to recover the 
solutions and gradients of the original and desired interface.  Note, the decisions
on when to use shifts and when to use flips are listed in \cite{shuAccurateGradientApproximation2014}.


\subsection{Our Proposed Method and Contributions}\label{ss:contributions}

\begin{figure}
\centering
	\includegraphics[width=0.7\textwidth]{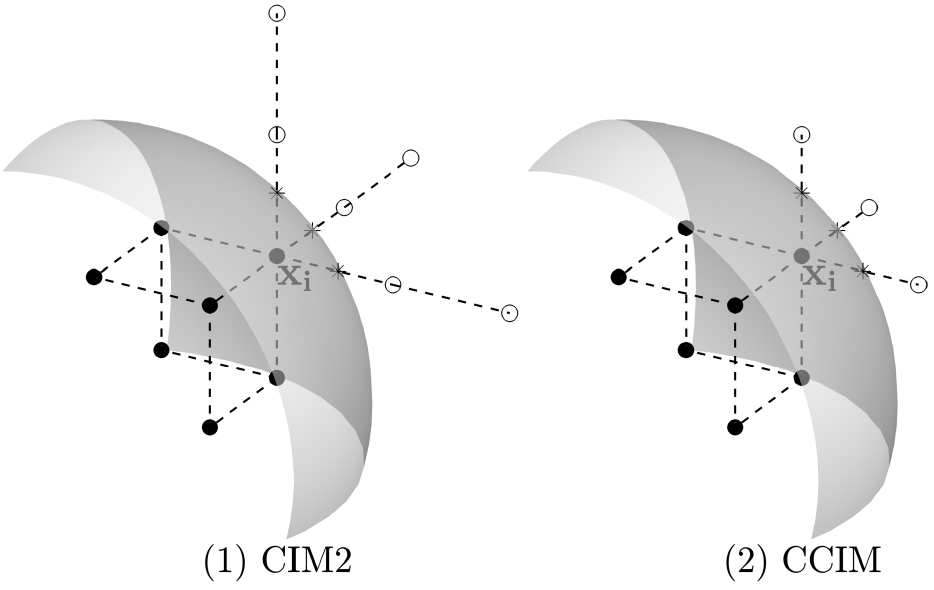}
\caption{Examples of a (1) CIM2 stencil  and a (2) CCIM stencil at $\vb{x_i}$. The circles and disks are grid points on different sides of the interface. The asterisks are the intersections of the surface and the grid lines. 
In this case, CIM2 and ICIM have the same stencil, which requires 2 points on the same side of the interface in each dimension \cite{chernCouplingInterfaceMethod2007}. Points that does not satisfy this requirement are handled in ICIM \cite{shuAccurateGradientApproximation2014}. The CCIM stencil requires fewer points so it's more compact. Both CIM2 and CCIM need some extra grid points compared to the standard central difference stencil.}
\label{f:2stencils}
\end{figure}

We propose what can be considered a hybrid that we call the Compact Coupling Interface 
Method (CCIM).  Our method combines elements of CIM \cite{chernCouplingInterfaceMethod2007} and 
Smereka's work on second-order accurate discrete delta functions by setting up an 
elliptic interface problem with interfacial jump conditions \cite{smerekaNumericalApproximationDelta2006}.
The use of Smereka's setup, itself based on Mayo's work in  \cite{mayoFastSolutionPoisson1984}, allows us to 
remove the quadratic polynomial approximations of CIM2 and its need for two points on 
either side of the interface in a direction that crosses 
the interface, thus compacting the stencil and allowing more applicability in generating 
accurate principal second-order derivatives.
Additional schemes are introduced to accurately handle mixed second-order derivatives by using more compact stencils on the 
same side or stencils from the opposite side, with the help of jump conditions, allowing
for the removal of exceptional points.  The result is that our constructed CCIM has a more compact stencil (shown in Fig~\ref{f:2stencils}) and it can approximate values and gradients of the solutions of elliptic interface problems with jump conditions with second-order accuracy in infinity norm for a variety 
interfaces, with observed advantages in robust convergence behavior in complex situations.

Note, while ICIM improves on CIM2 through shifting or flipping, the 
coupling equations remain largely the same, which only includes the principal second-order derivatives, and 
similar finite difference stencils \cite{shuAccurateGradientApproximation2014}. Our CCIM expands the 
coupling equations to include first-order derivatives as well, and utilizes more 
compact finite difference stencils, for the removal of additional exceptional points.

\subsection{Outline}

This paper is organized as follows. Section~\ref{ss:method} outlines the derivation and algorithm of CCIM. In Section~\ref{ss:result}, we show the convergence tests in three dimensions on geometric surfaces and two complex protein surfaces. We also test our method on a moving surface driven by the jump of the gradient at the interface. Section~\ref{ss:conclusion} is the conclusion.

\section{Method}\label{ss:method}

In $d$ dimensions, let $\Omega = [-1, 1]^d$ and discretize the domain uniformly with mesh size $h = 2/N$, where $N$ is the number of subintervals on one side of the region $\Omega$. Let $\vb{i} = (i_1,\ldots,i_d)$ be the multi-index with $i_k = 0,1,\ldots, N$ for $ k = 1,2,\ldots,d$. The grid points are denoted as $\bxi$ with the $k$-th coordinate $x_{k} = -1 + i_k h$. Let $\vb{e}_k$, $ k = 1,2,\ldots,d$ be the unit coordinate vectors. We also write $u(\vb{x_i}) = u_{\vb{i}}$. Here we use $\Delta u$ for the Laplacian of $u$ and $\hess u$ for the Hessian matrix of $u$. 
We use $\overline{ \vb{x_i} \vb x_{\vb i + \vb e_k}}$ to denote the grid segment between $\bxi$ and $\vb x_{\vb i + \vb e_k}$, and assume that the interface intersects with any grid segment at most once. 

Let $\bxi$ be a grid point at which we try to discretize the PDE. For notational simplicity, we drop the argument $\bxi$ and the dependency on $\vb i$ is implicit. We rewrite the PDE \eqref{eq:pde} at $\bxi$ as
\begin{equation}
	\label{eq:dimpde}
	- \sum_{k=1}^{d} \pdv{\epsilon}{x_k}  \pdv{u}{x_k} - \epsilon  \sum_{j=1}^{d} \pdv[2]{u}{x_k} + a u = f
\end{equation}

We classify the grid points into two categories:
if $\vb x_{\vb i - \vb e_k}$, $\bxi$ and $\vb x_{\vb i + \vb e_k}$ are in the same region in each coordinate direction, then we call $\bxi$ an \textbf{interior point}, otherwise $\bxi$ is called an \textbf{on-front point}. At interior points, standard central differencing gives a local truncation error of $\order{h^2}$ in  $\vb e_k$ direction. 
Our goal is to construct finite difference schemes with $\order{h}$ local truncation error at on-front points. The overall accuracy will still be second-order since the on-front points belong to a lower dimensional set \cite{levequeImmersedInterfaceMethod1994}. 
In this section, we derive a first-order approximation for the term $\pdv*{u}{x_k}$ and $\pdv*[2]{u}{x_k}$ in terms of u-values on neighboring grid points.
We denote the set of neighboring grid points of $\bxi$ as 
$B_r = \{\vb {x_j} \mid \norm{\vb j - \vb i}_\infty \leq r\}$
and call $r$ the radius of our finite difference stencil. 
As an example in \ref{f:2stencils}, the CIM2/ICIM stencil has $r = 2$ and the CCIM stencil has $r = 1$.

\subsection{Dimension-by-dimension discretization}\label{ss:dimbydim}

\begin{figure}
\centering
	\includegraphics[width=0.5\textwidth]{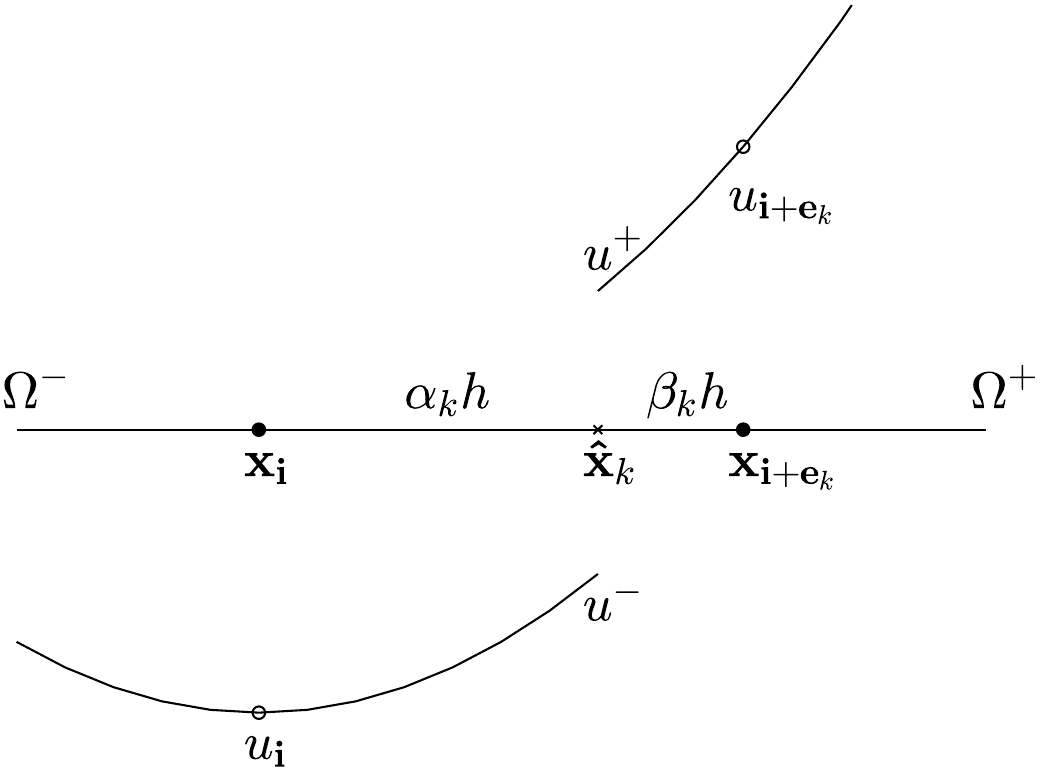}
\caption{The interface intersects the grid segment $\overline{ \vb{x_i} \vb x_{\vb i + \vb e_k}}$ at $\hxk$. $u^-$ and $u^+$ are the limits of u at $\hxk$ from $\omgm$ and $\omgp$. $u_{\vb i}$ and $u_{\vb{i}+\vb{e}_k}$ are approximated by Taylor's expansion at the interface.}
\label{f:diff}
\end{figure}

This section follows the derivation found in Smereka's work \cite{smerekaNumericalApproximationDelta2006}.
Along the coordinate direction $\vb{e}_k$, if the interface does not intersect the grid segment $\overline{ \vb{x_i} \vb x_{\vb i + \vb e_k}}$, then by Taylor's theorem,
\begin{equation}\label{eq:TaylorU}
u_{\vb i + \vb e_k} - u_{\vb i} 
		=  h\pdv{u}{x_k}   + \frac{h^2}{2}\pdv[2]{u}{x_k} +\order{h^3}.
\end{equation}
Suppose the interface intersects the grid segment $\overline{ \vb{x_i} \vb x_{\vb i + \vb e_k}}$ at the interface point $\hxk$. Let $\alpha_k = \norm{\hxk - \bxi}/h$ and $\beta_k = 1 - \alpha_k$. Suppose $\bxi$ is located in $\omgm$. Denote the limit of $u(\vb x)$ as $\vb x$ approaches $\hxk$ from $\omgm$ by $u^-$ , and the limit from the other side by $u^+$ (see Fig \ref{f:diff}).
By Taylor's theorem,
\begin{equation}
\begin{aligned}
u_{\vb i} & =  u^- - \alpha_k h \pdv{u^-}{x_k} + \frac{(\alpha_k h)^2}{2}\pdv[2]{u^-}{x_k} +\order{h^3},\\
u_{\vb i + \vb e_k} & =  u^+ + \beta_k h \pdv{u^+}{x_k} +\frac{(\beta_k h)^2}{2} \pdv[2]{u^+}{x_k} +\order{h^3}.\\
\end{aligned}  
\label{eq:TaylorInterface}
\end{equation}
Subtract the above two equations and write the right-hand side in terms of jumps and quantities from $\omgm$:
\begin{multline}\label{eq:OneSideInterface}
u_{\vb i + \vb e_k} - u_{\vb i} 
		=  [u] + \beta_k h\left[\pdv{u}{x_k}\right]+ h\pdv{u^-}{x_k}+\\
		\frac{h^2}{2}\beta_k^2 \left[\pdv[2]{u}{x_k}\right]+ 
			 \frac{h^2}{2}(\beta_k^2-\alpha_k^2)\pdv[2]{u^-}{x_k}+\order{h^3}. 
\end{multline}
We can approximate components of $\grad u^-$ and $\hess u^-$ by
\begin{equation}\label{eq:TaylorGradU}
\pdv{u^-}{x_i} = \pdv{u}{x_i} + \alpha_k h \pdv{u}{x_i}{x_k} + \order{h^2},
\end{equation}
\begin{equation}\label{eq:TaylorHessU}
\pdv{u^-}{x_i}{x_k} = \pdv{u}{x_i}{x_k} + \order{h}.
\end{equation}
with $1\leq i\leq k \leq d$.
Together with the given jump condition, $[u] = \tau$, \eqref{eq:OneSideInterface} can be written as 
\begin{equation}\label{eq:OneSideInterfaceSub}
\begin{aligned}
u_{\vb i + \vb e_k} - u_{\vb i} 
		=  \tau + \beta_k h\left[\pdv{u}{x_k}\right]
		+ h \left(\pdv{u}{x_k} + \alpha_k h \pdv[2]{u}{x_k}\right)\\
		+ \frac{h^2}{2}\beta_k^2 \left[\pdv[2]{u}{x_k}\right]+ 
			 \frac{h^2}{2}(\beta_k^2-\alpha_k^2)\pdv[2]{u}{x_k}
			  +\order{h^3}.
\end{aligned}
\end{equation}

We emphasize that all the first and second-order derivatives in \eqref{eq:OneSideInterfaceSub} are evaluated at the on-front point $\bxi$, and all the jump term $[\pdv*{u}{x_k}]$ and $[\pdv*[2]{u}{x_k}]$ are defined and evaluated at $\hxk$, which are intermediate quantities that will be eliminated in the coupling equation.

\subsection{Coupling Equation}\label{ss:couplingeqn}

This section then sets up coupling equations following the work of 
\cite{chernCouplingInterfaceMethod2007}.
In \eqref{eq:OneSideInterfaceSub}, suppose we can approximate the jump $[\pdv*{u}{x_k}]$ and $[\pdv*[2]{u}{x_k}]$ in terms of $u_{\vb j}$, $\pdv*{u}{x_k}$ and $\pdv*[2]{u}{x_k}$, with $1 \leq k \leq d$ and $\vb j \in B_r$ for some stencil radius $r$. 
Then in each coordinate direction, for $1\leq k \leq d$, we can write down two equations, \eqref{eq:OneSideInterfaceSub} or \eqref{eq:TaylorU}, by considering the two grid segments $\overline{ \vb{x_i} \vb x_{\vb i + s \vb e_k}}$ for $s = \pm 1$. In $d$ dimensions we have $2d$ equations and $2d$ unknowns: the first-order derivatives $\pdv*{u}{x_k}$ and the principal second-order derivatives $\pdv*[2]{u}{x_k}$ for $1\leq k \leq d$.  This leads to a system of linear equations of the following form:
\begin{equation}\label{eq:coupling}
  	M  \begin{pmatrix} \pdv{u}{x_k} \\ \pdv[2]{u}{x_k} \end{pmatrix}_{1\leq k\leq d} = \frac{1}{h^2}  \begin{pmatrix} L_{k,s}(\ujbr) \end{pmatrix}_{1\leq k \leq d, s = \pm 1} + 
  	\order{h}
  \end{equation}
where $L_{k,s}(\ujbr)$ is some affine function of $u$-values in the neighborhood $B_r$. We call \eqref{eq:coupling} the coupling equation and $M$ the coupling matrix. By inverting $M$, we can approximate $\pdv*{u}{x_k}$ and $\pdv*[2]{u}{x_k}$ in terms of $u$-values and obtain the finite difference approximation of the PDE \eqref{eq:dimpde} at the on-front point $\vb{x_i}$.

In the next few sections, we describe the ingredients to construct the coupling equation \eqref{eq:coupling}. 
We will use ``y = L(x)'' to denote “write the quantities y in terms of affine function of
quantities x”, where L represent a generic affine function and it depends on the geometric quantities of the interface.
In Section~\ref{ss:jumpdu}, we derive expressions to approximate the jump of the first-order derivatives $[\pdv*{u}{x_k}]$ in \eqref{eq:OneSideInterfaceSub}. In Section~\ref{ss:jumpd2u}, we approximate the jump of the principal second-order derivatives $[\pdv*[2]{u}{x_k}]$ in \eqref{eq:OneSideInterfaceSub}. In Section~\ref{ss:mixdu}, we discuss how to approximate the mixed derivatives, which is used to approximate $[\pdv*{u}{x_k}]$ and $[\pdv*[2]{u}{x_k}]$. In Section~\ref{ss:algo}, we combine all the ingredients and describe our algorithm to obtain the coupling equation.

\subsection{Approximation of $[\pdv*{u}{x_k}]$} \label{ss:jumpdu}
Let $\vb n$ be the unit normal vector at the interface, and $\vb s_1$, \ldots, $\vb s_{d-1}$ be the unit tangent vectors. 
The tangent vectors can be obtained by projecting the coordinate vectors onto the tangent plane.
We can write
\begin{equation}\label{eq:jumpgradu}
[\grad u] 
= [\grad u \cdot \vb n] \vb n + \sum_{j=1}^{d-1}[\grad u \cdot \vb s_j] \vb s_j
= [\grad u \cdot \vb n] \vb n + \sum_{j=1}^{d-1}(\grad \tau \cdot \vb s_j) \vb s_j,
\end{equation}
because $[\grad u \cdot \vb{s}_j] = \grad \tau \cdot \vb s_j$ for $1\leq j \leq d-1$. 

We use the following trick frequently to decouple the jump in subsequent derivations:
\begin{equation}\label{eq:trick}
\left[\epsilon v\right] 
		=  \epsilon^+[v]+[\epsilon]v^-.
\end{equation}
The jump condition $[\epsilon\nabla u\cdot \vn] = \sigma$ can be rewritten as
\begin{equation}\label{eq:jumpDun}
[\grad u \cdot \vb n] = \frac{1}{\epsilon^+} (\sigma - [\epsilon] \grad u^- \cdot \vb n).
\end{equation}
Substitute \eqref{eq:jumpDun} into \eqref{eq:jumpgradu}, in direction $\vb e_k$, we have
\begin{equation}\label{eq:jumpux}
\left[\pdv{u}{x_k}\right] = \frac{1}{\epsilon^+} (\sigma - [\epsilon] \grad u^- \cdot \vb n) (\vb n \cdot \vb e_k) + \sum_{j=1}^{d-1} (\grad \tau \cdot \vb s_j) (\vb s_j \cdot \vb e_k).
\end{equation}
Approximate $\grad u^-$ by Taylor's theorem \eqref{eq:TaylorGradU}, we get 
\begin{equation}\label{eq:jumpuxSub}
\begin{aligned}
\left[\pdv{u}{x_k}\right] = \frac{1}{\epsilon^+} \left(\sigma - [\epsilon] \sum_{j=1}^d \left( \pdv{u}{x_j} + \alpha_k h \pdv{u}{x_j}{x_k}\right) (\vb n \cdot \vb e_j)  \right) (\vb n \cdot \vb e_k) \\ + \sum_{j=1}^{d-1} (\grad \tau \cdot \vb s_j) (\vb s_j \cdot \vb e_k).
\end{aligned}
\end{equation}
Notice that the jump in the first-order derivative $[\pdv*{u}{x_k}]$ can be written as linear combinations of the first-order derivatives $\pdv*{u}{x_l}$, $1\leq l \leq d$ and the second-order derivatives $\pdv*{u}{x_l}{x_m}$, $ 1\leq l\leq m$:
\begin{equation}\label{eq:jumpduform}
\left[\pdv{u}{x_k}\right] = L( \grad u, \hess u).
\end{equation}
where $L$ represents the affine function in \eqref{eq:jumpuxSub}.
Our goal is to approximate the mixed derivatives $\pdv*{u}{x_k}{x_l}$, $k\neq l$ in terms of the neighboring $u$-values $u_{\vb j}$, $\vb j \in B_r$, the first-order derivatives $\pdv*{u}{x_k}$, and the principal second-order derivatives $\pdv*[2]{u}{x_k}$, $1\leq k \leq d$ (see Section ~\ref{ss:jumpd2u} and~\ref{ss:mixdu}), which are the terms used in the coupling equation \eqref{eq:coupling}.



\subsection{Approximation of $[\pdv*[2]{u}{x_k}]$} 
\label{ss:jumpd2u}
To remove the jump of the principal second-order derivatives $[\pdv*[2]{u}{x_k}]$, $ k = 1,2,\ldots,d$ in \eqref{eq:OneSideInterfaceSub}, we need to solve a system of linear equations, whose unknowns are all the jump of the principal and the mixed second-order derivatives. This idea is used in Smereka's work \cite{smerekaNumericalApproximationDelta2006} to arrive at the discrete approximation of the delta function. 
The detailed derivation is given in the Appendix~\ref{ss:derivation}. However, we would like to note here our definition of the gradient of a vector field as there is no standard one. 
For a vector field $\vb v$, we define the matrix $\grad\vb v$ to have entries in row $i$ and column $j$:
$$
(\grad \vb v)_{ij} = \frac{\partial v_j}{\partial x_i}.
$$
This particular form is used in our equations on the jumps of second derivatives.

\textbf{Tangential derivative of jump of tangential derivative:}
The first set of equations are obtained by differentiating the interface boundary condition in the tangential directions. For $m = 1,\cdots,d-1$ and $n = m,\cdots,d-1$, we get $d(d-1)/2$ equations
\begin{equation}
\grad [\grad u \cdot \vb s_m] \cdot \vb s_n = \grad (\grad \tau \cdot \vb s_m) \cdot \vb s_n.
\end{equation}
Expand each term, and with the help of \eqref{eq:trick} and \eqref{eq:jumpux}, we get
\begin{equation} \label{eq:dtau}
\vb{s}_n^T [\grad^2u]\vb{s}_m = 
\vb{s}_n^T \grad^2 \tau \vb{s}_m 
- \frac{1}{\epsilon^+}( \sigma - [\epsilon] \grad u^- \cdot \vb n) \vb{s}_n^T \grad \vn \vb{s}_m 
- (\grad \tau \cdot \vn ) \vb{s}_n^T \grad \vn \vb{s}_m.
\end{equation}

\textbf{Tangential derivative of flux jump:}
By differentiating the jump of flux in the tangential directions, we get another $d-1$ equations for $m = 1,\cdots,d-1$, 
\begin{equation}
\grad[\epsilon \grad u \cdot \vb n] \cdot \vb s_m = \grad \sigma \cdot \vb s_m  .
\end{equation}
After expansion,
\begin{equation}\label{eq:dsigma}
\begin{aligned}
\vb{s}_m^T [\grad^2u]\vn 
&= \frac{1}{\epsilon^+} \grad \sigma \cdot \vb{s}_m 
- \frac{[\epsilon]}{\epsilon^+} \vb{s}_m^T \hess u^- \vn 
- \frac{[\epsilon]}{\epsilon^+} \vb{s}_m^T \grad \vn \grad u^- \\
&- \sum_{k=1}^{d-1} (\grad\tau\cdot \vb{s}_k)  \vb{s}_m^T \grad \vn \vb{s}_k 
- \frac{1}{(\epsilon^+)^2} (\grad \epsilon^+ \cdot \vb{s}_m ) ( \sigma - [\epsilon] \grad u^- \cdot \vb n)\\
&- \frac{1}{\epsilon^+}  [\grad \epsilon \cdot \vb{s}_m]  (\grad u^- \cdot \vb n).
\end{aligned}
\end{equation}

\textbf{Jump of PDE:}
The final equation comes from the jump of the PDE:
\begin{equation}
	[-\grad \cdot (\epsilon \grad u) + au] = [f].
\end{equation} 
After expansion, we have
\begin{equation}\label{eq:df}
\begin{aligned}
[\Delta u] 
&= - \left[\frac{f}{\epsilon}\right] 
+ \frac{a^+}{\epsilon^+}\tau 
+ \left[\frac{a}{\epsilon}\right]u^-\\
&- \frac{1}{\epsilon^+} \sum_{k=1}^{d-1} (\grad\tau \cdot \vb{s}_k) (\grad \epsilon^+ \cdot \vb{s}_k)
- \left[\frac{\grad \epsilon}{\epsilon}\right] \cdot \grad u^-.
\end{aligned}
\end{equation}
Combining \eqref{eq:dsigma}, \eqref{eq:dtau}, and \eqref{eq:df}, we arrive at a system of linear equations whose unknowns are the jump of the second-order derivatives:
\begin{equation}\label{eq:jump}
G \left(\left[\pdv{u}{x_k}{x_l}\right]\right)_{1 \leq k \leq l \leq d} = L\left(u^-, \grad u^-, \hess u^- \right).
\end{equation}
where $G$ is a matrix that only depends on the normal and the tangent vectors, $L$ stands for the affine function in \eqref{eq:dtau}, \eqref{eq:dsigma} and \eqref{eq:df}.
In two and three dimensions, it can be shown that the absolute value of the determinant of $G$ is 1 (See Appendix~\ref{ss:derivation}).

As an example of \eqref{eq:jump} in two dimensions, let $\vb{s} = [s_1, s_2]^T$ and $\vn = [n_1, n_2]^T$, and assume that $\epsilon(\vb{x})$ is a piecewise constant function, then the system of linear equations \eqref{eq:jump} is given by
\begin{equation}\label{eq:jump2d}
\begin{aligned}
&
\begin{pmatrix}
s_1^2 & s_2^2 & 2s_1s_2\\
s_1n_1 & s_2n_2 & s_1n_2+s_2n_1\\
1 & 1 & 0\\
\end{pmatrix}
\begin{pmatrix}
[u_{xx}]\\
[u_{yy}]\\
[u_{xy}]\\
\end{pmatrix}
=\\
&\begin{pmatrix}
\vb{s}^T \hess \tau \vb{s} 
- \frac{1}{\epsilon^+} \left(\sigma-[\epsilon] \nabla u^-\cdot \vn \right) \vb{s}^T\grad \vn \vb{s} 
- (\grad \tau \cdot \vb n) \vb{s}^T \grad \vn \vb{s}\\
\frac{1}{\epsilon^+}\grad \sigma \cdot \vb{s} 
- \frac{[\epsilon]}{\epsilon^+}( \vb{s}^T \hess u^- \vn + \vb{s}^T \grad \vn \grad u^- ) 
- (\grad \tau \cdot \vb{s}) \vb{s}^T\grad \vn \vb{s}\\
-[\frac{f}{\epsilon}] 
+ \frac{a^+}{\epsilon^+} \tau 
- [\frac{a}{\epsilon}] u^-\\
\end{pmatrix}.
\end{aligned}	
\end{equation}

By Taylor's theorem as in \eqref{eq:TaylorU}, \eqref{eq:TaylorGradU} and \eqref{eq:TaylorHessU}, $u^-$, $\grad u^-$ and $\hess u^-$ can all be approximated by $u_{\vb j}$, $\vb j \in B_r$ and components of $\grad u$ and $\hess u$ at the grid point. Therefore, after substitution,  \eqref{eq:jump} has the form
\begin{equation}\label{eq:jumpd2uform}
G \left(\left[\pdv{u}{x_k}{x_l}\right]\right)_{1 \leq k \leq l \leq d} = 
L\left( \ujbr, \grad u, \hess u \right)
\end{equation}
where $ \vb{j} \in B_r$ and $L$ represent an updated affine function based on \eqref{eq:jump} after substitution.


Recall that our goal is to write every quantities in terms of $u_{\vb j}$, $\vb j \in B_r$, $\pdv*{u}{x_k}$, $\pdv*[2]{u}{x_k}$, $1\leq k\leq d$, called ``allowable terms'', which are the terms in the coupling equation \eqref{eq:coupling}. Here we lay out the steps to achieve this goal.
Firstly, we approximated the mixed derivatives $\pdv*{u}{x_k}{x_l}$, $k\neq l$, using the allowable terms and the jump of the mixed derivatives (see Section~\ref{ss:mixdu}):
\begin{equation} \label{eq:apprxmix}
\pdv{u}{x_k}{x_l} =  L \left( \ujbr, \pdv{u}{x_k}, \pdv[2]{u}{x_k}, \left[\pdv[2]{u}{x_k}{x_l}\right] \right)
\end{equation}
where $L$ represent the finite difference scheme to approximate the mixed derivatives.
Secondly, we substitute the above expressions \eqref{eq:apprxmix} into \eqref{eq:jumpd2uform} to eliminate the mixed derivatives:
\begin{equation}\label{eq:jumpd2uform2}
	G \left(\left[\pdv{u}{x_k}{x_l}\right]\right)_{1 \leq k \leq l \leq d} = 
	L\left(\ujbr, \grad u, \left(\pdv[2]{u}{x_k}\right)_{1\leq k\leq d} \right)
\end{equation}
with an updated affine function $G$ and $L$ matrix based on \eqref{eq:jumpd2uform}.
Notice that we have eliminated the mixed derivatives in $\hess u$ in \eqref{eq:jumpd2uform}.
Now the right hand side of \eqref{eq:jumpd2uform2} only contains the allowable terms. 
By solving the linear system of equations, we can write the jump of the second-order derivatives $[\pdv*[2]{u}{x_k}{x_l}]$, $1\leq k\leq l \leq d$ using the allowable terms.

If jump of the second-order derivatives is used in \eqref{eq:apprxmix}, then we can substitute the above expressions \eqref{eq:jumpd2uform2} into \eqref{eq:apprxmix} to eliminate the jump of the second-order derivatives, $[\pdv*[2]{u}{x_k}{x_l}]$, $1\leq k < l \leq d$.
This leads to
\begin{equation} \label{eq:apprxmix2}
	\pdv{u}{x_k}{x_l}=  
	L\left(\ujbr, \pdv{u}{x_k}, \pdv[2]{u}{x_k} \right)
\end{equation}
with some affine function $L$ obtained by eliminating the jump of the second-order derivatives in \eqref{eq:apprxmix}.
Now the right hand side of \eqref{eq:apprxmix2}, and hence \eqref{eq:jumpuxSub}, also only contains the allowable terms. 

Finally, we substitute the expressions of $[\pdv*{u}{x_k}]$, and $[\pdv*[2]{u}{x_k}]$, both only involve the allowable terms at this stage, into \eqref{eq:OneSideInterfaceSub}. After rearrangement, we obtain one row of the coupling equation \eqref{eq:coupling}.
The full algorithm will be summarized in Section~\ref{ss:algo}.

\subsection{Approximation of the mixed derivative}
\label{ss:mixdu}

\begin{figure}[!htbp]
\centering
\includegraphics[width=\textwidth]{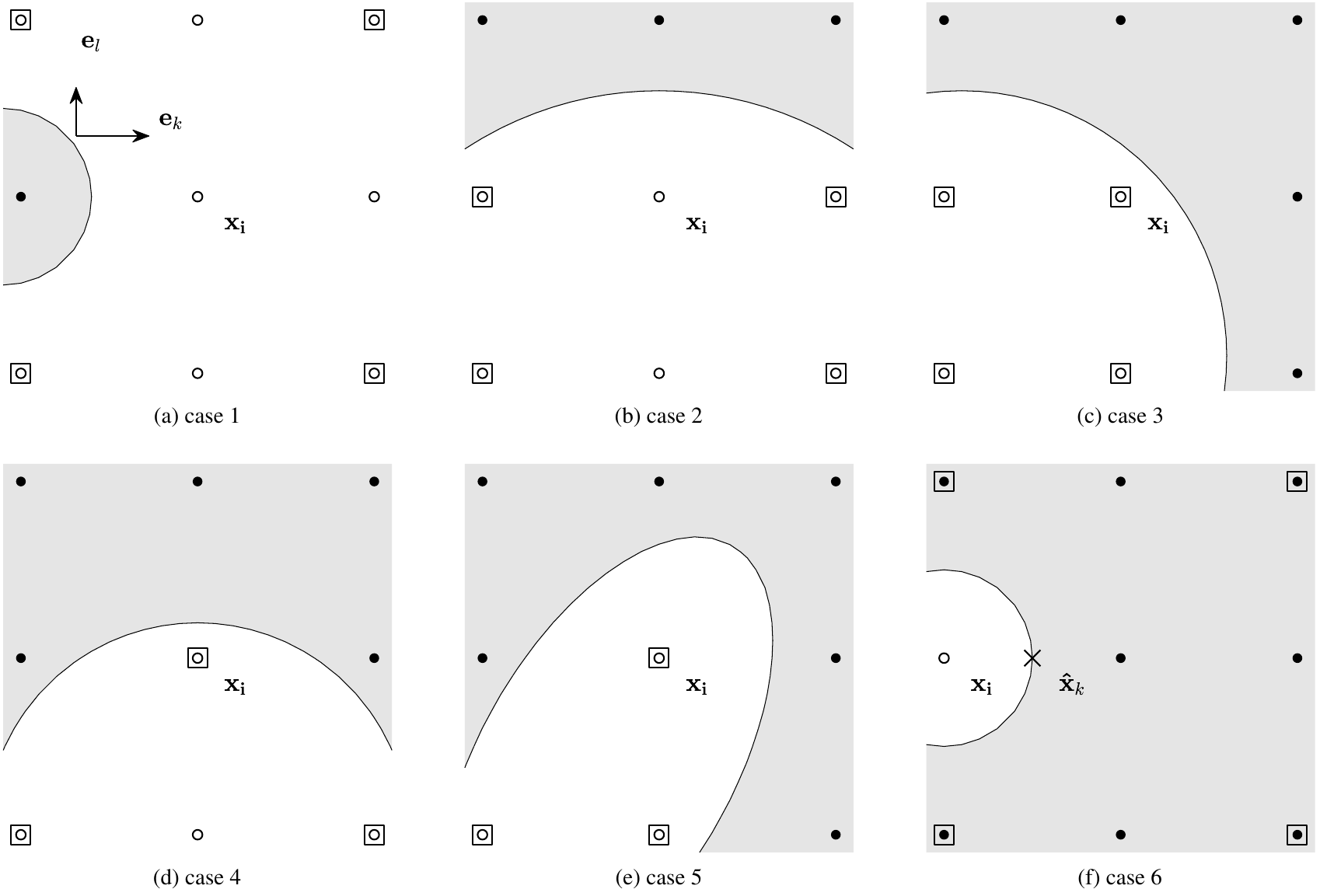}
\caption{Approximation of the mixed derivative $\pdv*[2]{u}{x_k}{x_l}$ at $\bxi$. The circles and disks are grid points in $\omgm$ and $\omgp$. The $u$-values at the squares are used to approximate the mixed derivative. Case 1 is the usual central difference. Case 2 and 3 are biased difference. Case 4 uses the first-order derivatives at at $\bxi$. Case 5 uses the first-order and the second-order derivatives at at $\bxi$. Case 6 uses the jump $[\pdv*[2]{u}{x_k}{x_l}]$ at $\hxk$ and the mixed second-order derivative on the other side at $\vb x_{\vb i + \vb e_k}$, which is approximated by central difference.}
\label{f:mixdu}
\end{figure}

Depending how the interface intersects the grid, different schemes are needed to approximate the mixed derivative $\pdv*{u}{x_k}{x_l}$, $k\neq l$ at $\bxi$. Notice that we are allowed to make use of the first-order and the second-order derivatives, as they are the variables in the coupling equation \eqref{eq:coupling}. 

In Section \ref{ss:mixducases}, we first describe the available finite difference schemes for the mixed derivatives in different scenarios. In Section \ref{ss:mixduorder}, we describe the decision to choose different schemes when multiple schemes are available. 

\subsubsection{Schemes for the mixed derivatives}
\label{ss:mixducases}

Though any $\order{h}$ approximation suffices, we prefer schemes with smaller local truncation error. Therefore, in all the following formula, we also compute the $\order{h}$ term explicitly. 
In Fig.~\ref{f:mixdu}, we demonstrate examples of different scenarios, and we describe the schemes in the following:

\textbf{Case 1 (central difference)}
\begin{equation}
\pdv{u}{x_k}{x_l} = \frac{1}{4h^2} \left( u_{\vb{i} + \vb{e}_k + \vb{e}_l} 
- u_{\vb{i} - \vb{e}_k + \vb{e}_l}
- u_{\vb{i} + \vb{e}_k - \vb{e}_l}
+ u_{\vb{i} - \vb{e}_k - \vb{e}_l} \right)+\order{h^2}.
\end{equation}

\textbf{Case 2 (biased differencing with rectangular stencil)}
\begin{equation}
\label{eq:uxyrect}
\begin{aligned}
 \pdv{u}{x_k}{x_l} &= \frac{1}{2h^2} \left( u_{\vb{i}+\vb{e}_k} - u_{\vb{i}+\vb{e}_k-\vb{e}_l} - u_{\vb{i}-\vb{e}_k} + u_{\vb{i}-\vb{e}_k-\vb{e}_l} \right)  \\
 & + \frac{1}{2} h \pdv[3]{u}{x_k}{x_l^2} + \order{h^2}
\end{aligned}
\end{equation}


\textbf{Case 3 (biased differencing with square stencil)}
\begin{equation}
\label{eq:uxysquare}
\begin{aligned}
 \pdv{u}{x_k}{x_l} & = \frac{1}{h^2} (u_{\vb{i}} - u_{\vb{i}-\vb{e}_k} - u_{\vb{i}-\vb{e}_l} + u_{\vb{i}-\vb{e}_k-\vb{e}_l})\\
 & + \frac{1}{2} h \left(\pdv[3]{u}{x_k}{x_l^2} + \pdv[3]{u}{x_k^2}{x_l}\right) + \order{h^2}
\end{aligned}
\end{equation}

\textbf{Case 4 (triangular stencil with first-order derivatives)}
In case 4, we can make use of the first-order derivatives
\begin{equation}
\label{eq:uxy1/4tri}
\begin{aligned}
\pdv{u}{x_k}{x_l} &= \frac{1}{2h^2} \left( 2h\pdv{u}{x_k} + u_{\vb{i}+\vb{e}_k-\vb{e}_l} - u_{\vb{i}-\vb{e}_k-\vb{e}_l}\right)\\
	& + \frac{1}{2} h \left(\pdv[3]{u}{x_k}{x_l^2} + \frac{1}{3} \pdv[3]{u}{x_k}\right) + \order{h^2}
\end{aligned}
\end{equation}

\textbf{Case 5 (triangular stencil with first and second-order derivatives)}
In case 5, we can make use of the first-order and the principal second-order derivatives
\begin{equation}
\label{eq:uxy1/8tri}
\begin{aligned}
\pdv{u}{x_k}{x_l} &= \frac{1}{h^2} \left( h \pdv{u}{x_k} - \frac{h^2}{2} \pdv[2]{u}{x_k} - u_{\vb{i}-\vb{e}_l} + u_{\vb{i}-\vb{e}_k-\vb{e}_l} \right) \\
& + \frac{1}{2} h \left(\pdv[3]{u}{x_k}{x_l^2} + \pdv[3]{u}{x_k^2}{x_l} + \frac{1}{3} \pdv[3]{u}{x_k}\right) + \order{h^2}
\end{aligned}
\end{equation}

\textbf{Case 6 (shifting to the other side)}
When there are not enough grid points on the same side, we can make use of the mixed derivative on the other side of the interface and the jump of the mixed derivative

\begin{equation}\label{eq:mixjump}
\begin{aligned}
\pdv{u}{x_k}{x_l} 
&= \pdv{u}{x_k}{x_l}\/(\vb{x}_{\vb i + \vb{e}_k}) - \left[\pdv{u}{x_k}{x_l}\right]_\hxk+\order{h}\\
\end{aligned}
\end{equation}
where $\pdv*{u}{x_k}{x_l}\/(\vb{x_i}+h\vb{e}_k)$ can be approximated by $u$-values on the other side of the interface using central difference as in case 1 to 3. 
In case 6, the finite difference stencil will have a radius $r = 2$, as $u$-values of more than one grid point away are used. 
For example, as shown in Fig.~\ref{f:mixdu} (case 6), $u(\vb{x}_{\vb i + 2 \vb{e}_k + \vb{e}_l})$ is used to approximate the mixed derivative at $\vb{x_i}$.

\textbf{Case 7 (shifting to the same side)}
Though not illustrated in Fig.~\ref{f:mixdu}, it's also possible to approximate the mixed derivative at $\vb{x_i}$ by the mixed derivative of $\vb{x_i}$'s direct neighbor that is on the same side of the region:
\begin{equation}\label{eq:shift}
\pdv{u}{x_k}{x_l} = \pdv{u}{x_k}{x_l}\/(\vb{x}_{\vb i + s \vb{e}_m}) + \order{h},
\end{equation}
where $s=\pm 1$ and $1\leq m\leq d$, and $\vb{x}_{\vb i + s \vb{e}_m}$ is on the same side of the interface as $\vb{x_i}$.
Then the $\pdv*{u}{x_k}{x_l}\/(\vb{x}_{\vb i + s \vb{e}_m})$ can be approximated using case 1, 2 or 3, which only involves $u$-values on the same side of the interface.
This is the same as the ``shifting'' strategy used in ICIM \cite{shuAccurateGradientApproximation2014}. If $m=k$ or $l$, that is, we are shifting in the $kl$-plane, then the stencil have a radius $r = 2$. Otherwise, we are shifting out of the $kl$-plane, and the stencil would have a radius $r = 1$. Therefore, shifting out of plane is preferred for a more compact stencil.

We note that Case 1, 2 and 3, which approximate the mixed derivatives by neighboring u-values, are the default choices in CIM and hence ICIM.
When such approximations cannot be found, Case 7 (shifting to the same side) is considered in ICIM.
Case 4 and 5 are unique to our CCIM, as they make use the first-order derivatives.
These two techniques are not allowed in the CIM and ICIM, because their coupling equations do not involve the first-order derivatives.
Case 6 (shifting to the other side) is also unique to our CCIM, as it makes use of the jump of the mixed derivatives.
These techniques (case 4, 5 and 6) help to make our stencil more compact, and allow us to handle very complicated interface even when the grid is relatively coarse.


\begin{figure}[!htbp]
\centering
\includegraphics[width=0.7\textwidth]{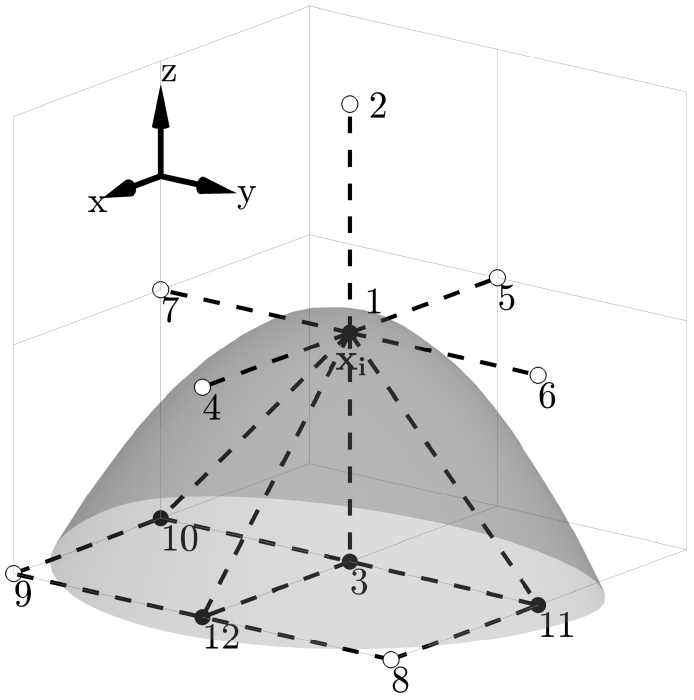}
\caption{An example of CCIM stencil at $\vb{x_i}$. The grid points in the stencil are labelled from 1 to 12. The disks and the circles are the grid points outside and inside the surface. (1,10,11) are used to approximate $u_{yz}(\vb{x_i})$ (see Fig.~\ref{f:mixdu} case 4). (1, 3, 12) are used to approximate $u_{zx}(\vb{x_i})$ (Fig.~\ref{f:mixdu} case 5). (8, 9, 10, 11) can be used to approximate $u_{xy}(\vb{x_i})$. Alternatively, $u_{xy}(\vb{x_i})$ can be approximated by using $u_{xy}$ at the other side (any one of 2, 4, 5, 6 or 7) and the jump condition (see Fig.~\ref{f:mixdu} case 6).}
\label{f:ccimexample}
\end{figure}



\subsubsection{Ordering of Different Schemes}
\label{ss:mixduorder}
Multiple schemes to approximate the mixed derivatives might be available at the same grid point, and a natural question is which one to choose. Here we discuss the criteria for choosing the scheme. Overall, we would like the scheme to be simple, compact and accurate, both in terms of the local truncation error and the condition number of the coupling matrix.

\textbf{Simplicity}
For simplicity, we prefer schemes that only use $u$-values as they are easy to implement. Therefore case 1, 2 and 3, which using central differencing or biased differencing are preferred. 
When the derivatives have to be involved, case 4 only use first derivatives, and is preferred over case 5, which also uses second derivatives.

\textbf{Compactness}
For compactness, we want to have a smaller radius $r$ for our finite difference stencil. Therefore case 6a (shifting to the other side) is the least preferable. When case 7 (shifting to the same side) has to be used, we prefer shifting out-of-plane than shifting in-plane.

\textbf{Local truncation error}
For accuracy, we look at the $\order{h}$ of the local truncation error. Central differencing (case 1) is preferred over biased differencing (and case 2 over case 3). Case 3 or case 4 has similar local truncation error.
In general, case 6 (shifting) will leads to larger local truncation error compared with case 1 to 5.

\textbf{Condition number of the coupling equation}
Another consideration for accuracy is the condition number of the coupling equation. Solving a linear system with large condition number is prone to large numerical errors. 
Therefore, in cases where both case 3 and case 4 are available, we choose the scheme that leads to the coupling matrix with a smaller estimated condition number computed by \cite{hagerConditionEstimates1984}.
The approximate condition number provide a good estimate of the actual condition number, which is computationally expensive to compute.
The effect is shown in Section~\ref{ss:result}. 

\textbf{Ordering}
As a summary of the above discussion on the criteria of choosing the differencing schemes for approximating the mixed derivatives,
here is how we rank the schemes from the most preferable to the least preferable: 
case 1 (central differencing);
case 2 (biased differencing with rectangular stencil);
case 3 (biased differencing with square stencil);
or
case 4 (triangular stencil with first-order derivatives), whichever leads to coupling matrix with smaller estimated condition number;
case 5 (triangular stencil with first and second-order derivatives);
case 7 (shifting to the same side);
case 6 (shifting to the other side); 
When multiple schemes are available to approximate the mixed derivatives, we use the most preferable scheme available. 
Though we can construct surfaces for a specific grid size such that none of the above schemes works, for smooth surfaces we can refine the grid such that the above schemes suffice. We note that case 5 and case 6 can be removed by refining the grid, while case 4 cannot, which is proved in \cite{shuAccurateGradientApproximation2014}.

\subsection{Algorithm}\label{ss:algo}
We describe our method to obtain the coupling equation at an on-front point in algorithmic order in Algorithm~\ref{algo}. Once we have the coupling equations \eqref{eq:coupling}, by inverting the coupling matrix, $\pdv*{u}{x_k}$ and $\pdv*[2]{u}{x_k}$, $1\leq k \leq d$ can be approximated by linear functions of $u_{\vb j}$, $\vb j \in B_r$

\begin{algorithm}[!htbp]
	\caption{Coupling equation at an on-front point $\bxi$}
	\label{algo}
	\textbf{Notation:} we use $y = L(x)$ to denote ``write the quantities y in terms of affine function of quantities x'', where $L$ represent a generic affine function.
	\begin{algorithmic}[1]
		\FOR{$1 \leq k\leq d$}
			\FOR{$s = \pm 1$}
			\IF {the interface intersects $\overline{ \vb{x_i} \vb x_{\vb i + s \vb e_k}}$ at $\hxk$}
	
			\FOR{$1 \leq j \leq d$, $j \neq k$}
			\STATE
				As described in Section~\ref{ss:mixdu}
				\[\pdv{u}{x_k}{x_j} = L\left(\ujbr, \pdv{u}{x_k}, \pdv[2]{u}{x_k}, \left[\pdv{u}{x_k}{x_j}\right]\right)\]
				
			\ENDFOR
	
			\STATE As described in Section~\ref{ss:jumpd2u},
			\[ \left[\pdv{u}{x_k}{x_j}\right]_{1\leq k \leq j\leq d} = 
			L\left(\ujbr, \left(\pdv{u}{x_j}\right)_{1\leq j\leq d}, \left(\pdv[2]{u}{x_j}\right)_{1\leq j\leq d}\right)
			\]
			
			\STATE By back substitution
			\[
				\left(\pdv{u}{x_k}{x_j}\right)_{1\leq j\leq d} =
				L\left(\ujbr, \left(\pdv{u}{x_j}\right)_{1\leq j\leq d}, \left(\pdv[2]{u}{x_j}\right)_{1\leq j\leq d}\right)
			\]

			\STATE Substitute the expression for the mixed derivatives into \eqref{eq:jumpuxSub}
			\[
				\left[\pdv{u}{x_k}\right] = L\left(\ujbr, \left(\pdv{u}{x_j}\right)_{1\leq j\leq d}, \left(\pdv[2]{u}{x_j}\right)_{1\leq j\leq d}\right)
			\]

			\STATE
			Substituting $[\pdv*{u}{x_k}]$ and $[\pdv*[2]{u}{x_k}]$ into \eqref{eq:OneSideInterfaceSub}, and after rearrangement, this gives one row of the coupling equations \eqref{eq:coupling}.

			\ELSE
			\STATE Direct application of Taylor's theorem \eqref{eq:TaylorU} gives one row of the coupling equations \eqref{eq:coupling}.
			\ENDIF
			\ENDFOR
		\ENDFOR
	\end{algorithmic}
\end{algorithm}

To get more stable convergence results, at grid points where case 1 and 2 not available, but case 3 and 4 are available, we use the algorithm to obtain two systems of coupling equations, and choose the system with a smaller estimated condition number of the coupling matrix. The effect of this criterion is demonstrated in Section~\ref{ss:eg1}. 

Let $A \vb{u}=f$ be the system of linear equations obtained from the finite difference approximation of the PDE \eqref{eq:dimpde}, where $\vb u$ is the vector of unknown $u$-values at the grid points. 
The algorithm to assemble the system of linear equations $A \vb{u}=f$ is summarized in Algorithm~\ref{algo:all}. 

\begin{algorithm}[!htbp]
	\caption{Linear systems of equations for the elliptic interface problem}
	\label{algo:all}
	\begin{algorithmic}[1]
		\FOR{Every grid point $\bxi$}
			\IF {$\bxi$ is an on-front point}
				\STATE Use Algorithm~\ref{algo} to obtain the coupling equation at $\bxi$. Solve the coupling equation to obtain the finite difference approximations of $\pdv*{u}{x_k}$ and $\pdv*[2]{u}{x_k}$, $1\leq k \leq d$.
			\ELSE
				\STATE Use standard central difference to approximate $\pdv*{u}{x_k}$ and $\pdv*[2]{u}{x_k}$, $1\leq k \leq d$.
			\ENDIF

			\STATE 
			Substitute the approximations of $\pdv*{u}{x_k}$ and $\pdv*[2]{u}{x_k}$,  $1\leq k \leq d$, into the PDE \eqref{eq:dimpde}, and obtain one row of the linear system of equations.
		\ENDFOR
	\end{algorithmic}
\end{algorithm}

\section{Numerical results}\label{ss:result}

We test our method in three dimensions with different surfaces. The first set of tests contains six geometric surfaces that are used in \cite{shuAccurateGradientApproximation2014}. And the second set of tests uses two complex biomolecular surfaces. These two sets are compared with our implementation of ICIM \cite{shuAccurateGradientApproximation2014} with the same setup. As tests in \cite{shuAccurateGradientApproximation2014} do not include the $a(\vb x)$ term, the third set of tests are the same six geometric surfaces with the $a(\vb x)$ term. The last test is a sphere expanding under a normal velocity given by the derivative of the solution in normal direction. Let $u_e$ be the exact solution of \eqref{eq:pde}, and $u$ be the numerical solution. 

For tests with a static interface, we look at the maximum error of the solution at all grid points, denoted as $\norm{u_e-u}_\infty$, and the maximum error of the gradient at all the intersections of the interface and the grid lines, denoted as $\norm{\grad u_e - \grad u}_{\infty,\Gamma}$. For the expanding sphere, we look at the maximum error and the Root Mean Square Error (RMSE) of the radius at all the intersections of the interface and the grid lines. All the tests are performed on a 2017 iMac with 3.5 GHz Intel Core i5 and 16GB memory. We use the AMG method implemented in the HYPRE library \cite{falgoutHypreLibraryHigh2002} to solve the sparse linear systems to a tolerance of $10^{-9}$.

In our work, the interface is represented by a level set function.
To obtain the location of the interface on the grid segment, we use a degree-6 interpolating polynomial of the level set function and find the root using regula-falsi method.
For our numerical tests, and in many applications \cite{shuAccurateGradientApproximation2014,zhouVariationalImplicitSolvation2014,zhongImplicitBoundaryIntegral2018}, we have the expression of $\tau(x,y,z)$ and $\sigma(x,y,z)$.
The quantities $\nabla \tau$, $\nabla^2 \tau$, $\nabla \sigma$, the normal vectors and the tangent vectors are all approximated numerically by central differencing at the grid point, and then interpolated to the interface location.
We believe that our method only requires these geometric quantities to be approximated to second-order accuracy for computational purposes. And this is verified empirically in our numerical tests.
We note that since we assume interface $\Gamma$ and the interfacial jump condition $\tau$ and $\sigma$ are all smooth,
we can consider any smooth extension of $\tau$ and $\sigma$ off the interface. 
Even if such smooth extensions are not available in an analytical form, the extensions can be computed to second order accuracy by algorithm such as the fast marching method \cite{sethianFastMarchingLevel1996}.

\subsection{Example 1}
\label{ss:eg1}
We test several geometric interfaces as in \cite{shuAccurateGradientApproximation2014}. The surfaces are shown in Fig.\ref{f:6surf}. Their level set functions are given below:
\begin{itemize}
	\item Eight balls: $\phi(x,y,z) = \min_{0\leq k\leq 7} \sqrt{(x-x_k)^2+(y-y_k)^2+(z-z_k)^2} - 0.3$, where $(x_k,y_k,z_k) = ((-1)^{\floor{k/4}}\times 0.5,(-1)^{\floor{k/2}}\times 0.5,(-1)^k\times 0.5)$
	\item Ellipsoid: $\phi(x,y,z)=2x^2 +3y^2 +6z^2-1.3$
	\item Peanut: $\phi(x,y,z) = \phi(r,\theta,\psi) = r - 0.5 - 0.2\sin(2\theta)\sin(\psi)$
	\item Donut: $\phi(x,y,z) = (\sqrt{x^2+y^2}-0.6)^2+z^2-0.4^2$
	\item Banana: $\phi(x,y,z) = (7x+6)^4+2401y^4+3601.5z^4+98(7x+6)^2(y^2+z^2)+4802y^2z^2-94(7x+6)^2+3822y^2-4606z^2+1521$
	\item Popcorn: 
	\begin{align*}
	\phi(x,y,z) 
	&= \sqrt{x^2 +y^2 +z^2} -r_0 \\
	& -\sum_{k=0}^{11} \exp(25((x-x_k)^2 +(y-y_k)^2 +(z-z_k)^2))
	\end{align*}
	where 
	\begin{align*}
	&(x_k,y_k,z_k) \\
	&= \frac{r_0}{\sqrt{5}}
	\left( 2\cos \left(\frac{2k\pi}{5}-\floor{\frac{k}{5}}\right),
	2\sin \left(\frac{2k\pi}{5}-\floor{\frac{k}{5}} \right),
	(-1)^{\floor{\frac{k}{5}}}\right), 0\leq k\leq 9\\
	&=r_0(0,0,(-1)^{k-10}),10\leq k\leq 11.
	\end{align*}
\end{itemize}
The computational domain is $\Omega=[-1,1]^3$. 
Let $\omgp = \{\phi<0\}$ and $\omgm = \{\phi>0\}$, which are the interior and exterior of the surface respectively.
The boundary condition is $u = 0$ on the boundary of the computational domain. 
The setup in this example is the same as that in \cite{chernCouplingInterfaceMethod2007,shuAccurateGradientApproximation2014}: we take $a=0$ and use use the same exact solution \eqref{eq:testu} and the coefficient \eqref{eq:epsilon}:
\begin{equation}
\label{eq:testu}
	u_e(x,y,z) = 
		\begin{cases}
			xy+x^4+y^4+xz^2+\cos(2x+y^2+z^3) & \text{if $(x,y,z)\in\omgp$}\\
			x^3+xy^2+y^3+z^4+\sin(3(x^2+y^2)) & \text{if $(x,y,z)\in\omgm$}
		\end{cases}       
\end{equation}
and
\begin{equation}
\label{eq:epsilon}
	\epsilon(x,y,z) = 
		\begin{cases}
			\epsp & \text{if $(x,y,z)\in\omgp$}\\
			\epsm & \text{if $(x,y,z)\in\omgm$}
		\end{cases}.
\end{equation}
with $\epsp = 80$ and $\epsm = 2$.
The source term and the jump conditions are calculated accordingly.

\begin{figure}[!htbp]
\centering
\includegraphics[width=\textwidth]{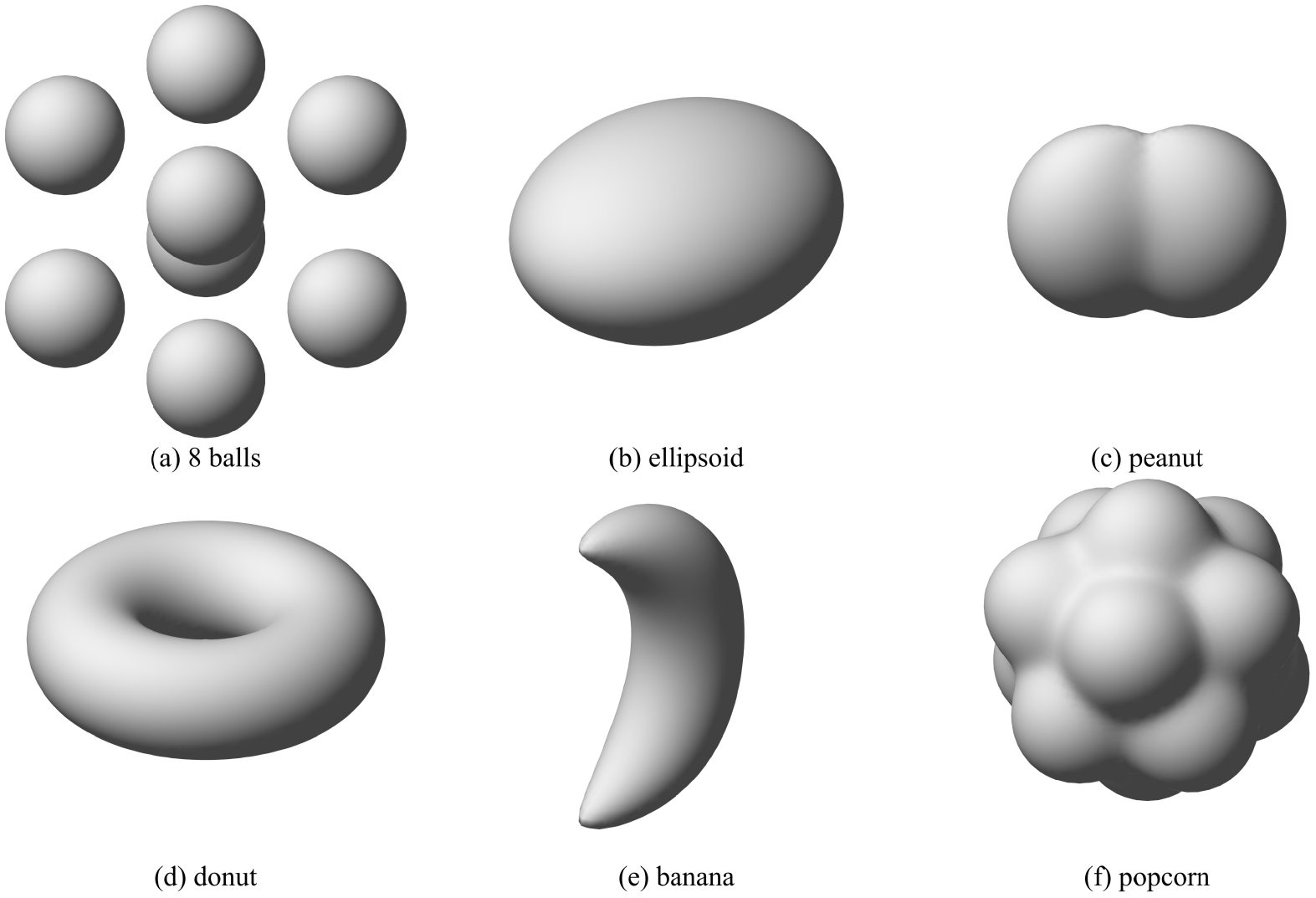}
\caption{The six interfaces (a) eight balls; (b) ellipsoid; (c) peanut; (d) donut; (e) banana; (f) popcorn}
\label{f:6surf}
\end{figure}

\begin{figure}[!htbp]
\centering
\includegraphics[width=\textwidth]{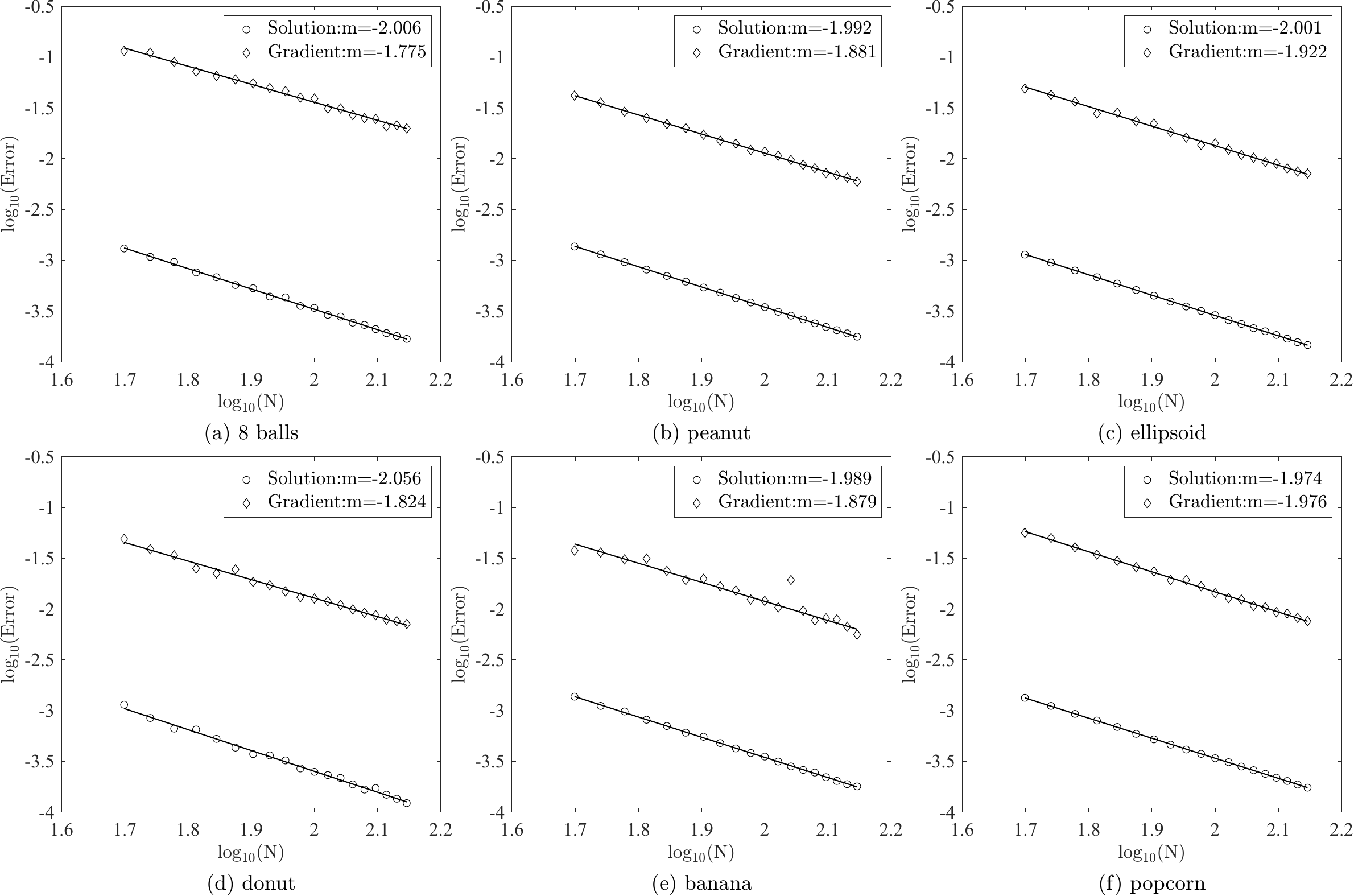}
\caption{The log-log plot of the error versus $N$ for the six surfaces. In each figure, N ranges from 50 to 140 with the increment $\Delta N = 5$. Circles are the maximum errors of the solution $\norm{u_e-u}_\infty$. Diamonds are the maximum errors of the gradient at interface $\norm{\grad u_e - \grad u}_{\infty,\Gamma}$. $m$ is the slope of the fitting line.}
\label{f:6surfconv}
\end{figure}

Fig \ref{f:6surfconv} shows the convergence result of the six interfaces. The convergence of the solution at grid points is second-order, and the convergence of the gradient at the interface is close to second-order. 

\subsubsection{Effect of different schemes for the mixed derivatives}
\label{ss:effect}
Next we demonstrate the effect of choosing the approximation schemes for the mixed derivatives based on the estimated condition numbers of the coupling matrices. As mentioned in Section~\ref{ss:mixdu}, when both case 3 and case 4 are available to approximate the mixed derivatives, we choose the scheme that gives a smaller estimated condition number of the coupling matrix. We denote this scheme as ``CCIM''. Alternatively, we can fix the order of preference for different methods. In ``scheme 1'', we always prefer case 4  (triangular stencil with first-order derivatives) to case 3 (biased differencing with square stencil).

\begin{figure}[!htbp]
\captionsetup[sub]{font=footnotesize,labelfont=footnotesize}
		 \centering
		 \begin{subfigure}[b]{0.5\textwidth}
				\centering
				 \includegraphics[width=\textwidth]{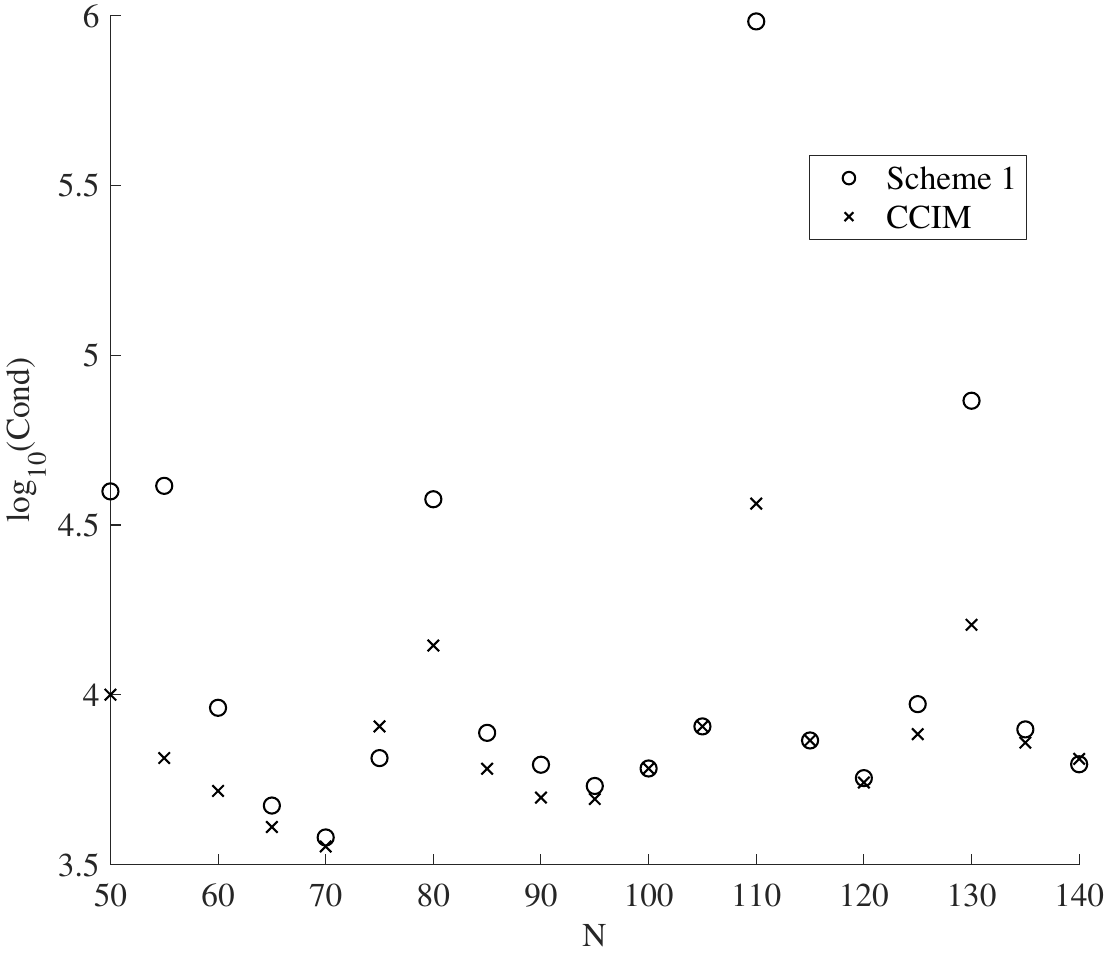}
				 \caption{maximum condition number}
				 \label{f:cond}
		 \end{subfigure}%
		 \hfill
		 \begin{subfigure}[b]{0.5\textwidth}
				\centering
				 \includegraphics[width=\textwidth]{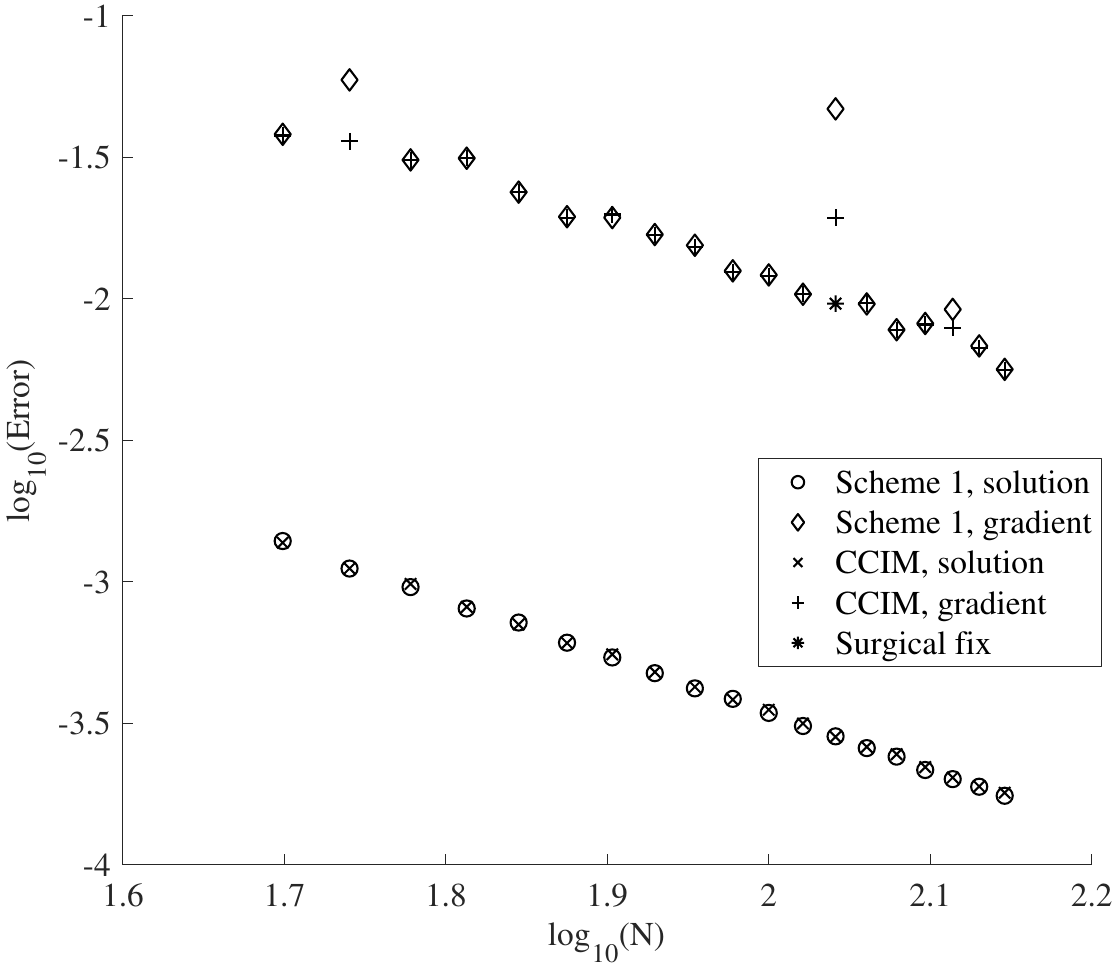}
				 \caption{Convergence}
				 \label{f:bananaconv}
		 \end{subfigure}
				\caption{Comparison of maximum condition numbers and convergence results with banana surface between scheme 1 and CCIM. In scheme 1, case 4 is preferred to case 3. In CCIM, case 3 and case 4 are chosen based on estimated condition numbers of the coupling matrices. (a) The maximum condition number of coupling matrices with different methods. (b) Log-log plot of the maximum errors in solution and gradient at the interface. 
				The relatively large error at $N=110$ for CCIM can be reduced by a surgical fix that choose the stencil with a smaller local truncation error.
				}
				\label{f:condandconv}
\end{figure}

Fig.~\ref{f:condandconv} demonstrates the effect of this decision using the banana shape surface as an example. For different $N$ and for both schemes, Fig.~\ref{f:cond} plots the maximum condition numbers (not the estimation from \cite{hagerConditionEstimates1984}) of all the coupling matrices, and Fig.~\ref{f:bananaconv} plots the convergence results of these two schemes. 
From Fig.~\ref{f:cond}, we can see that the maximum condition numbers in CCIM are almost always smaller than those in scheme 1 (except at $N = 75$, due to the estimation error \cite{hagerConditionEstimates1984}).
As shown in Fig.~\ref{f:bananaconv}, for most of the tests, CCIM and scheme 1 have roughly the same maximum error. We noticed that for $N = 110$, with scheme 1, at the interface point with the maximum error in the gradient, the coupling matrix has an exceptionally large condition number. 
By choosing the method with a smaller estimated condition number, we can get smaller error and obtain more stable convergence result in the gradient. 
If we prefer case 3 to case 4, then the results are similar to scheme 1: at some grid points large condition number is correlated to large error, and CCIM has more stable convergence behavior.

\subsubsection{Investigation of the outlier error}
\label{ss:outlier}
\begin{figure}[!htbp]
\centering
\includegraphics[width=0.5\textwidth]{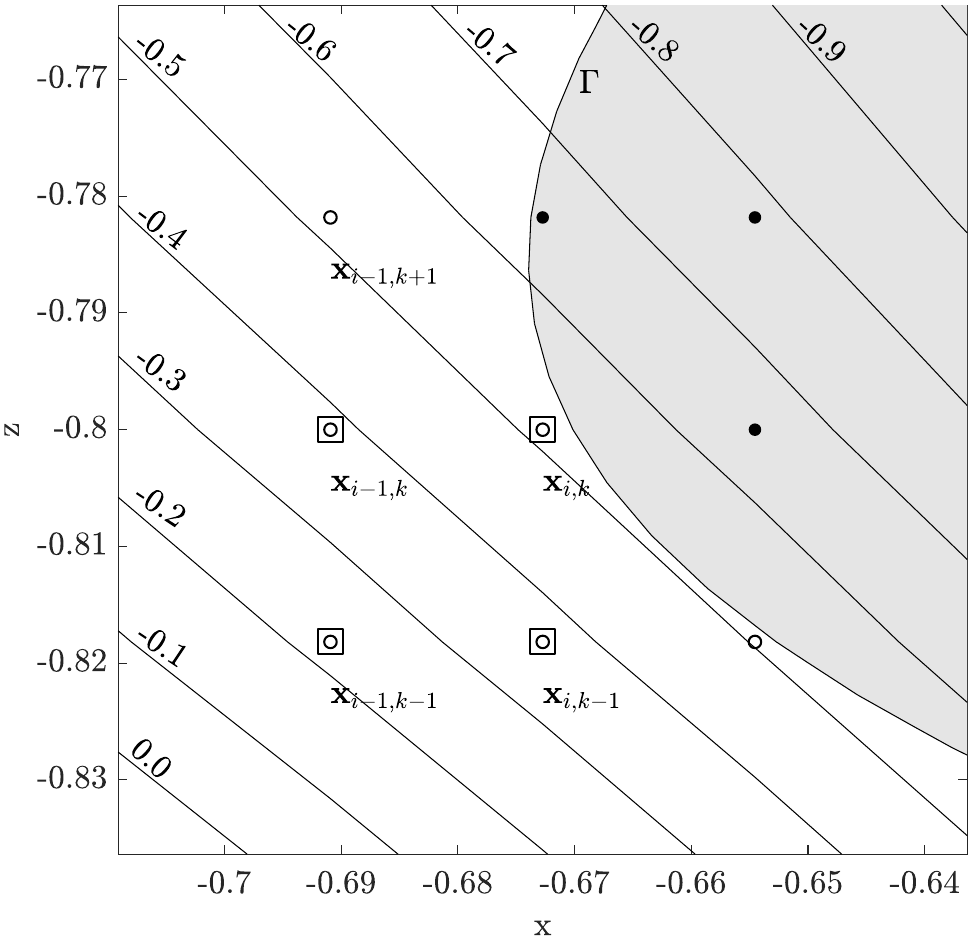}
\caption{Contour lines of $u_{xz}$ at the grid point where the maximum error in the gradient occurs for the banana interface at $N=110$ 
The disks and the circles are grid points outside and inside the surface. Grid points marked with squares (Fig.~\ref{f:mixdu} case 3) are used to approximate $u_{xz}(\vb{x}_{i,k})$ due to its simplicity, but has a local truncation error 0.160723,
The three points stencil (Fig.~\ref{f:mixdu} case 4) with $\vb{x}_{i,k}$, $\vb{x}_{i-1,k}$, $\vb{x}_{i-1,k+1}$ has a smaller local truncation error of 0.041853.
Use this stencil to approximate $u_{xz}(\vb{x}_{i,k})$ would reduce the error in the gradient at this grid point.
}
\label{f:jumpexplain}
\end{figure}

Though we can get a more stable convergence behavior by considering the condition number of the coupling matrices, there is an outlier of the error of the gradient for the banana interface at $N=110$ in Fig.~\ref{f:bananaconv}. A detailed analysis of the error reveals that it is caused by relatively large local truncation error when approximating $u_{xz}$. Fig.~\ref{f:jumpexplain} shows the contour line of the mixed derivative $u_{xz}$. 
However, due to the alignment of the surface with the grid, at $\vb{x}_{i,k}$, our algorithm uses the 4-point stencil $\vb{x}_{i,k}$, $\vb{x}_{i-1,k}$, $\vb{x}_{i,k-1}$ and $\vb{x}_{i-1,k-1}$ to approximate $u_{xz}(\vb{x}_{i,k})$ and has a local truncation error 0.160723. If we use the three point stencil $\vb{x}_{i,k}$, $\vb{x}_{i-1,k}$, $\vb{x}_{i-1,k+1}$, the local truncation error would be 0.041853, and the coupling matrix does not have large condition number. 
With this surgical fix, the final error would be in line with the rest of the data points, as shown in Fig~\ref{f:bananaconv} at N = 110, marked as ``Surgical fix''. This type of outliers happens rarely and does not affect the overall order of convergence. We apply this surgical fix only at this specific grid point to demonstrate a possible source of large error in the gradient.

In summary, though the overall order of convergence is second-order no matter which scheme is used to approximate the mixed second-order derivatives, a relatively large error can be caused by a large condition number of the coupling matrix, or a large local truncation error when approximating the mixed second-order derivative. When different schemes to approximate the mixed second-order derivatives are available, ideally we prefer the scheme that produces smaller local truncation error and smaller condition number of the coupling matrix. However, these two goals might be incompatible sometimes. It's time consuming to search through all the available schemes and find the one that leads the smallest condition number of the coupling matrix. It's also difficult to tell a priori which scheme gives smaller local truncation error. Therefore we try to find a middle ground by only considering the condition number when both case 1 and 2 are not available but case 3 and case 4 are available.

\begin{figure}[!htbp]
\centering
\includegraphics[width=\textwidth]{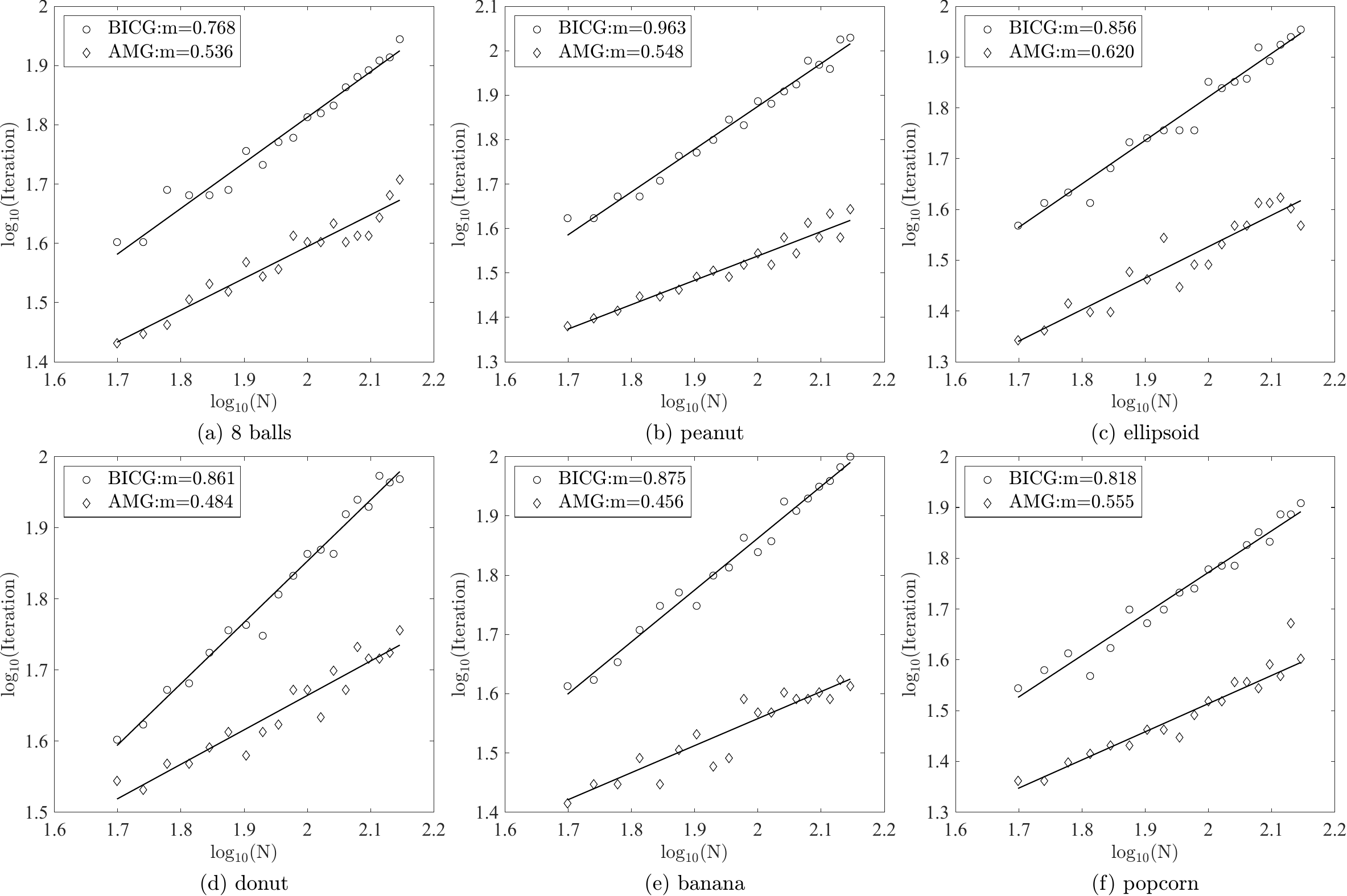}
\caption{The log-log plot of the number of iterations versus $N$ for the six surfaces. In each figure, N ranges from 50 to 140 with the increment $\Delta N = 5$. Circles are the numbers of iterations using Biconjugate Gradient Stabilized Method with ILU preconditioner (abbreviated as BICG). Diamonds are the numbers of iterations using AMG. $m$ is the slope of the fitting line. Though AMG has better scaling (smaller $m$), we report that both solvers take approximately the same CPU time in this range of $N$. }
\label{f:iteration}
\end{figure}

\begin{figure}[!htbp]
	\centering
	\includegraphics[width=\textwidth]{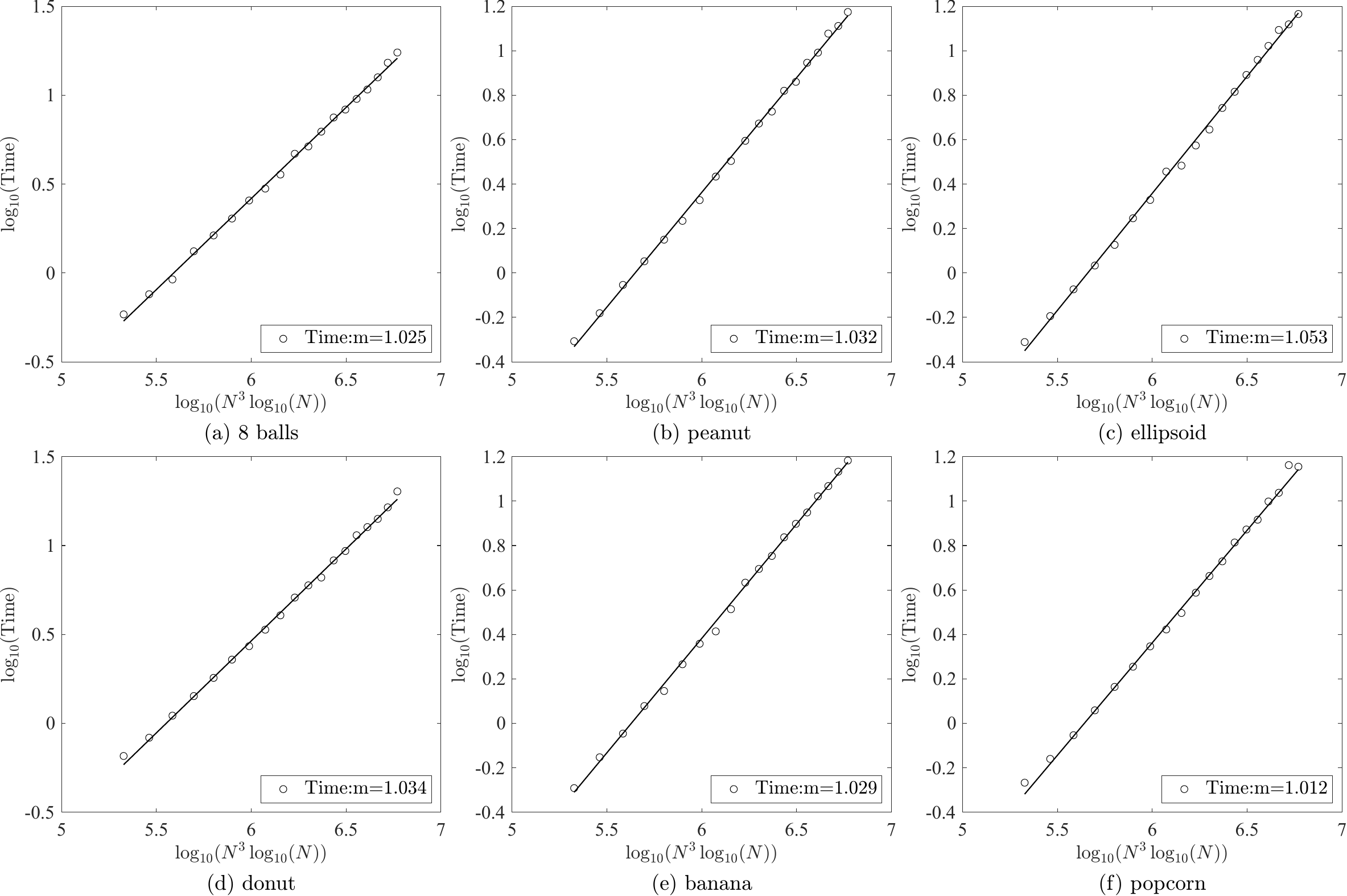}
	\caption{The log-log plot of the CPU time versus $N^3\log_{10}(N)$ for the six surfaces using AMG. $m$ is the slope of the fitting line. }
	\label{f:cpu}
\end{figure}

The resulting linear system for the PDE is sparse and asymmetric, and can be solved with any ``black-box'' linear solvers. Fig.~\ref{f:iteration} shows the log-log plot for the number of iterations versus $N$. We used Biconjugate Gradient Stabilized Method with ILU preconditioner (abbreviated as BICG) and Algebraic Multigrid Method (AMG), both are implemented in the HYPRE library \cite{falgoutHypreLibraryHigh2002}. The number of iterations grows linearly with $N$ for BiCGSTAB and sub-linearly for AMG. Though AMG has better scaling property, for the range of $N$ in Fig \ref{f:iteration}, both solvers take approximately the same CPU time. 
Fig.~\ref{f:cpu} shows the log-log plot of the CPU time versus $N^3\log_{10}(N)$ for the six surfaces using AMG. The slopes of the regression lines are close to 1. This scaling behavior is on par with that in ICIM (see Fig.14 in \cite{shuAccurateGradientApproximation2014}). All the tests are performed on an Apple M2 Max CPU.

Empirically, our method does not encounter issue when the interface points are very close to the grid points: the minimum $\alpha$ among all the test cases in Fig. \ref{f:6surfconv} is 2.3282e-18. 
In addition, in Appendix~\ref{ss:highcontrast}, we show that our method can handle high contrast problems with a large difference in the coefficients $\epsp$ and $\epsm$.

\subsection{Example 2}
\label{ss:eg2}

Next we test our method on two complex molecular surfaces and compare CCIM with our implementation of ICIM \cite{shuAccurateGradientApproximation2014}. The solvent accessible surface describes the interface between solute and solvent. Such interfaces are complex and important in applications. We construct the surfaces as in \cite{shuAccurateGradientApproximation2014}: from the PDB file of 1D63 \cite{brownCrystalStructureBerenild1992} and MDM2 \cite{kussieStructureMDM2Oncoprotein1996}, we use the PDB2PQR \cite{dolinskyPDB2PQRAutomatedPipeline2004} software to assign charges and radii using the AMBER force field. The PQR files contain information of the positions $\vb{p_i}$ and radii $r_i$ of the atoms. We scale the positions and radii such that the protein fit into our computation box. Then we construct the level set function of the interface as the union of smoothed bumps:
\begin{equation}
	\phi(\vb x) = c - \sum_{i} \chi_{\eta}(r_i - \norm{\vb x-\vb{p_i}}),
\end{equation}
where $\chi_{\eta}$ is a smoothed characteristic function
\begin{equation}
	\chi_{\eta} = \frac{1}{2}\left( 1+\tanh \left(\frac{x}{\eta}\right)\right).
\end{equation}

\begin{figure}[!htbp]
\centering
\includegraphics[width=\textwidth]{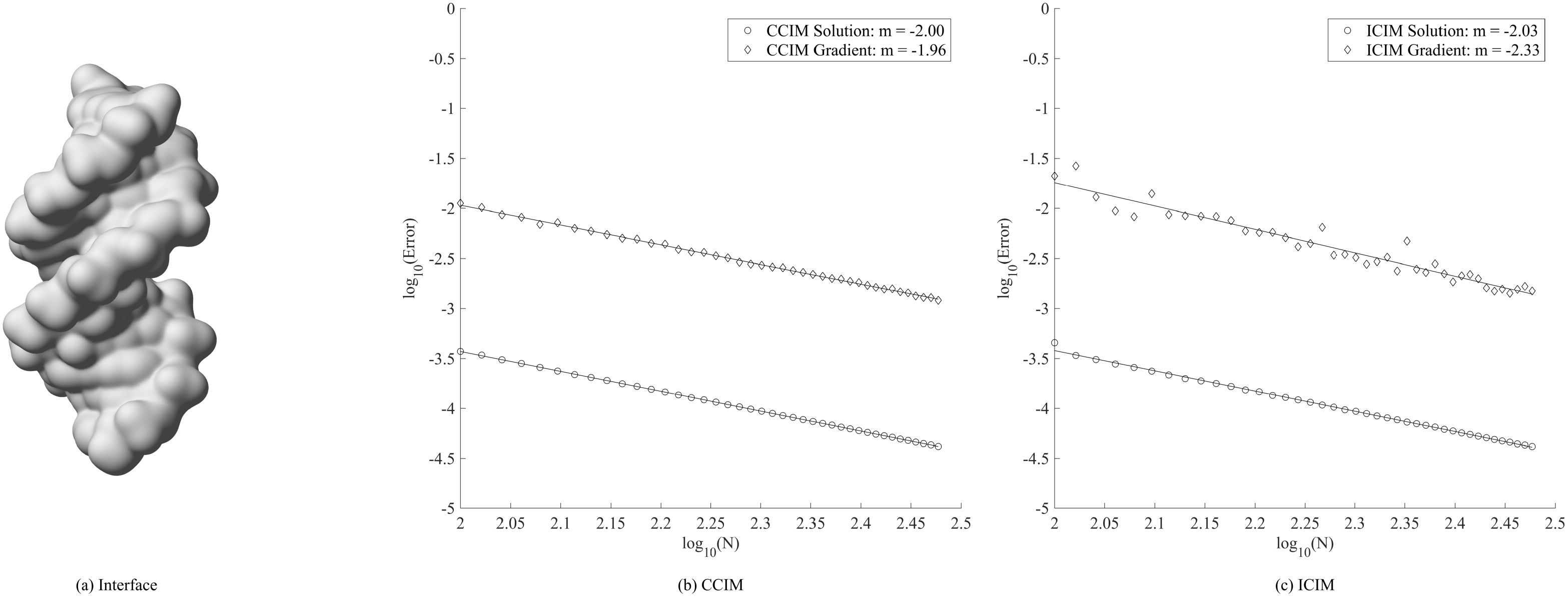}
\caption{Convergence result for 1D63 interface with $c=0.25$ and $\eta=1/40$. (a) The smooth surface of 1D63. (b) log-log plot of error by CCIM. (b) log-log plot of error by ICIM. $N$ ranges from 100 to 340 with the increment $\Delta N = 5$.}
\label{f:1d63eta40}
\end{figure}

\begin{figure}[!htbp]
\centering
\includegraphics[width=\textwidth]{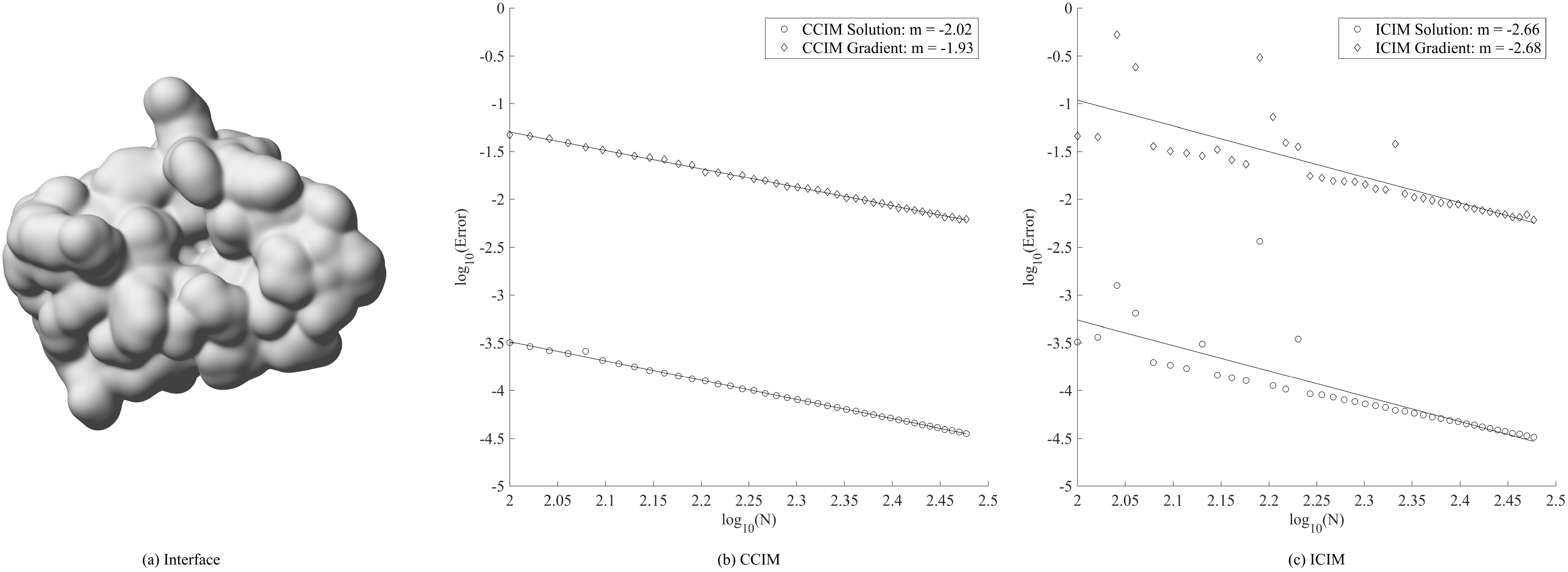}
\caption{Convergence result for MDM2 interface with $c=0.25$ and $\eta=1/30$. (a) The smooth surface of MDM2. (b) log-log plot of error by CCIM. (b) log-log plot of error by ICIM. $N$ ranges from 100 to 340 with the increment $\Delta N = 5$.}
\label{f:mdm2}
\end{figure}

The molecule 1D63 has 486 atoms and has a double-helix shape, as shown in Fig.~\ref{f:1d63eta40}(a). MDM2 has 1448 atoms, and the surface has a deep pocket to which other proteins can bind, as shown in Fig \ref{f:mdm2}(a). We also implement ICIM \cite{shuAccurateGradientApproximation2014} and compare the convergence results between CCIM and ICIM in Fig.~\ref{f:1d63eta40} and Fig.~\ref{f:mdm2}.

As shown in Fig.~\ref{f:1d63eta40} and Fig.~\ref{f:mdm2}, compared with our implementation of ICIM, the convergence results of CCIM is very robust even for complex interfaces. There is little fluctuation in the convergence results. In our ICIM implementation, the order of convergence exceeds second-order because large errors at coarse grid points skew the fitting line to have a more negative slope. The results demonstrate the advantage of the compactness in our CCIM formulation when dealing with complex surfaces.

\subsection{Example 3}
\label{ss:eg3}

\begin{figure}[!htbp]
\centering
\includegraphics[width=\textwidth]{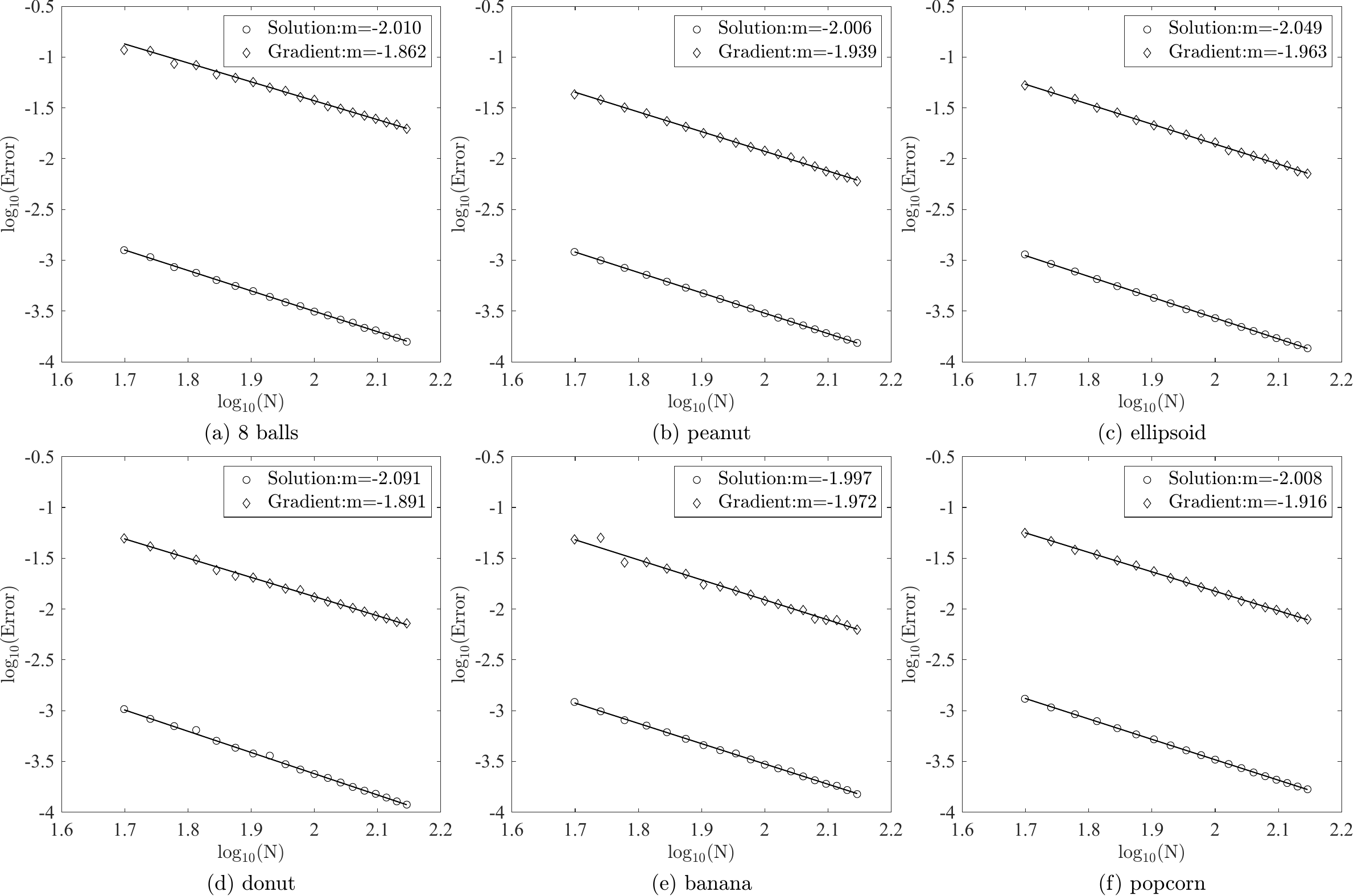}
\caption{The log-log plot of error with $a(\vb x)$ term versus $N$ for the six surfaces. 
In each figure, N ranges from 50 to 140 with the increment $\Delta N = 5$. 
The grids are perturbed by some uniformly distributed random numbers between 0 and mesh size $h$ in each coordinate.
Circles are the maximum errors of the solution $\norm{u_e-u}_\infty$. 
Diamonds are the maximum errors of the gradient at interface $\norm{\grad u_e - \grad u}_{\infty,\Gamma}$. $m$ is the slope of the fitting line.
}
\label{f:6surfconvwitha}
\end{figure}

We also test our problem with the same exact solution $u_e$ \eqref{eq:testu} and coefficients $\epsilon$ \eqref{eq:epsilon}, but with a non-zero $a(x,y,z)$ term:
\begin{equation}
	a(x,y,z) = 
		\begin{cases}
			2  \left(\sin(x) + e^{y^2}\right) & \text{if $(x,y,z)\in\omgm$}\\
			80 \left(\cos(z) + e^{-y^2}\right) & \text{if $(x,y,z)\in\omgp$}.
		\end{cases}      
\end{equation}
This term is not tested in CIM \cite{chernCouplingInterfaceMethod2007} and ICIM\cite{shuAccurateGradientApproximation2014}. 
In Fig.~\ref{f:6surfconvwitha}, we show the convergence result of the six surfaces.
The grids are perturbed by some uniformly distributed random numbers between 0 and mesh size $h$ in each coordinate.
This changes the alignment of the surface with the grid lines compared with Fig.~\ref{f:6surfconv}. 
We can see the convergence of the solution at grid points and the convergence of the gradient at the interface are second-order.

\subsection{Example 4}
\label{ss:eg4}

In this example we look at the evolution of an interface driven by the jump of the normal derivative of the solution using the level set method \cite{osherLevelSetMethods2003a}. Suppose the surface $\Gamma$ is evolved with normal velocity $v_n = [\grad u \cdot n]$. $\phi = \phi(\vb{x},t)$ is a level set function representing the evolving surface $\Gamma = \Gamma(t)$, i.e., $\Gamma(t) = \{\vb x\mid \phi(\vb x, t) = 0\}$. The dynamics of the interface is given by the level set equation,
\begin{equation}
	\label{eq:levelseteq}
	\phi_t + v_n |\grad \phi| = 0.
\end{equation}
We use the forward Euler method for first-order accurate time discretization, Godunov scheme for the Hamiltonian, and the Fast Marching Method \cite{sethianFastMarchingLevel1996} to extend $v_n$ to the whole computational domain.








We start with the radially symmetric exact solution
\begin{equation}
	u_e(\vb x) =
	\begin{cases}
	\frac{1}{1+\norm{\vb x}^2} & \text{$\vb x\in \omgm$} \\
	-\frac{1}{1+\norm{\vb x}^2} & \text{$\vb x\in\omgp$} \\
	\end{cases}
\end{equation}
and
\begin{equation}
	a(\vb x) =
	\begin{cases}
	2\sin(\norm{\vb x}) & \text{$\vb x\in \omgm$} \\
	80\cos(\norm{\vb x}) & \text{$\vb x\in \omgp$} \\
	\end{cases}.
\end{equation}
The coefficient $\epsilon$ is the same as \eqref{eq:epsilon}. The source term and the jump conditions are calculated accordingly.

If the surface is a sphere of radius $r$, by symmetry, the normal velocity is uniform over the sphere and is given by
\begin{equation}
	v_n(r) = [\grad{u}\cdot n] = \frac{4r}{(1+r^2)^2}.
\end{equation}
Let the initial surface be a sphere of radius 0.5, then the motion of the surface is described by the ODE
\begin{equation}
\frac{dr}{dt} = v_n(r),\quad r(0) = 0.5
\end{equation}
which can be computed to high accuracy. The result is a sphere expanding at varying time-dependent speeding.

\begin{figure}[!htbp]
\captionsetup[sub]{font=footnotesize,labelfont=footnotesize}
\centering
 \begin{subfigure}[b]{0.5\textwidth}
		\centering
		 \includegraphics[width=\textwidth]{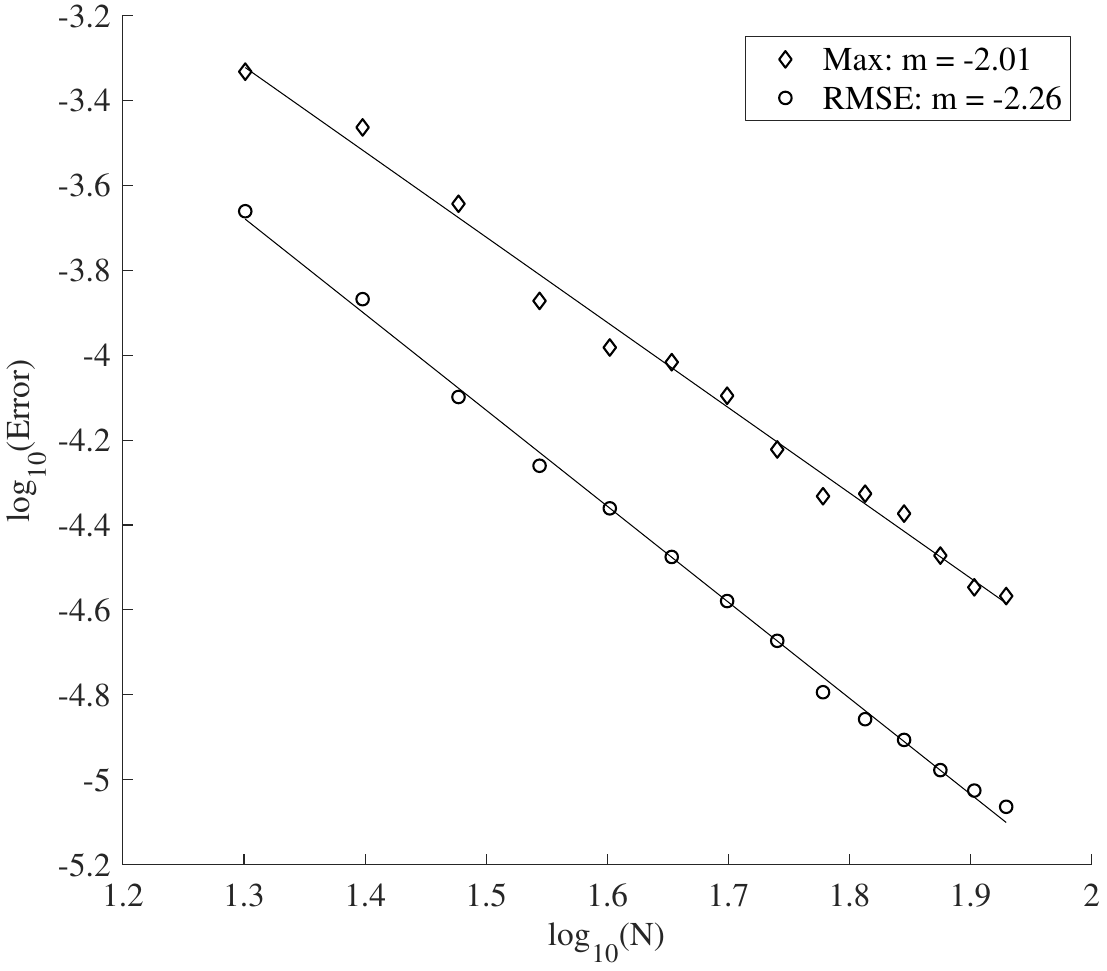}
		 \caption{Error for radii at final time.}
		 \label{f:motionerr}
 \end{subfigure}%
 \hfill
 \begin{subfigure}[b]{0.5\textwidth}
		\centering
		 \includegraphics[width=\textwidth]{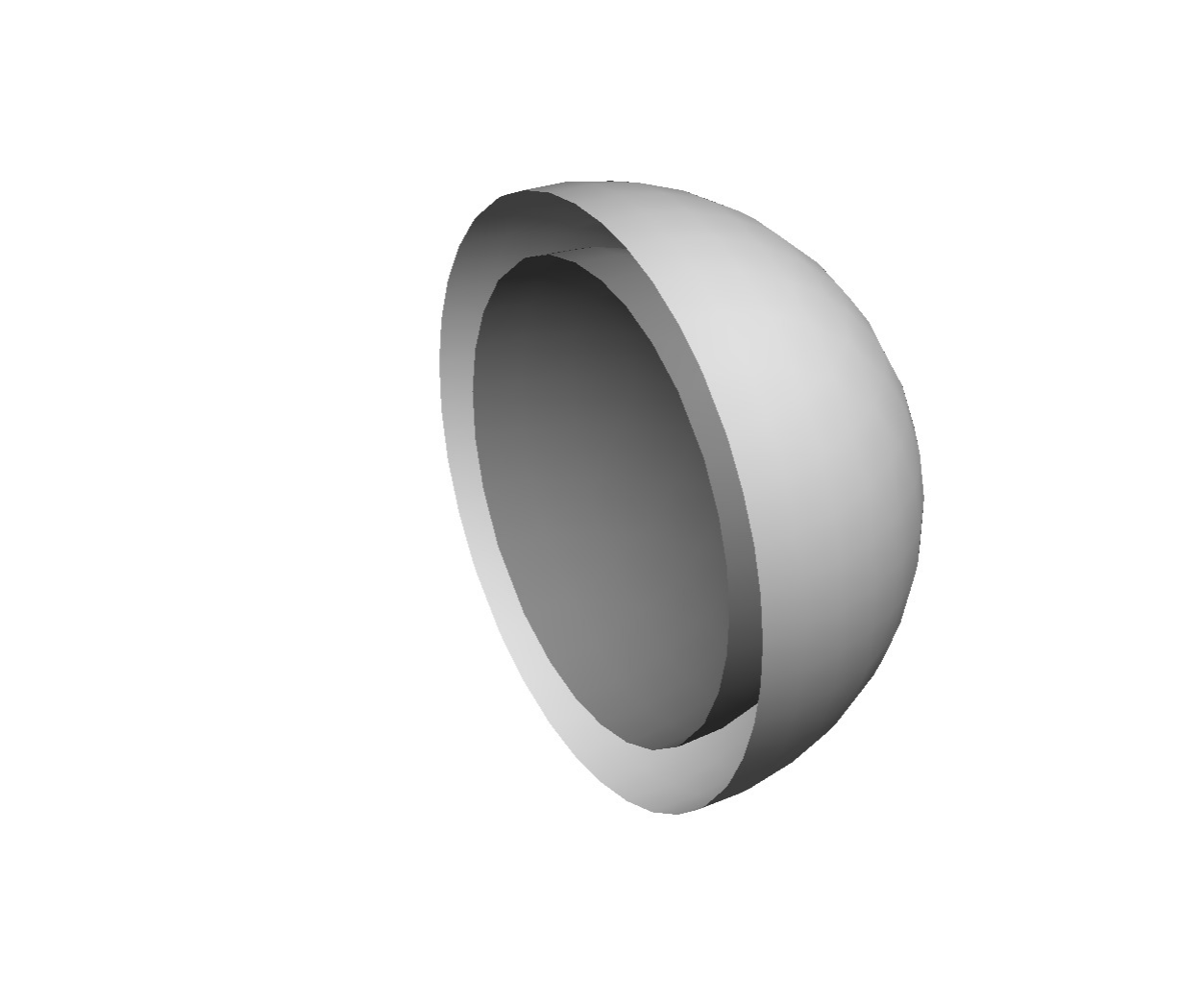}
		 \caption{Initial (inner) and final (outer) surface}
		 \label{f:motionsurf}
 \end{subfigure}
		\caption{(a) Maximum error and Root Mean Squared error (RMSE) of the radii measured at all the intersections of the surface and the grid lines. $N$ ranges from 20 to 85 with the increment $\Delta N = 5$. (b) Initial surface (inner) and final surface (outer) for $N$ = 20.}
		\label{f:motion}
\end{figure}

In Fig.~\ref{f:motionerr}, we look at the maximum error and the Root Mean Squared Error (RMSE) of all the radii obtained from the intersections of the surface and the grid lines at the final time $t = 0.1$ for different grid size $N$. 
The results are second-order accurate. In Fig.~\ref{f:motionsurf}, we plot the initial and final surface for $N = 20$. The shape is well-preserved. Without accurate gradient approximation, the surface might become distorted or oscillatory.

\section{Conclusions}\label{ss:conclusion}

In this paper, we proposed the Compact Coupling Interface Method (CCIM) to solve elliptic interface boundary value problems in any dimension. Our method combines elements from the Coupling Interface Method (CIM) and Mayo's approach to Poisson's equation on irregular regions. Standard central difference schemes are used at interior points. At on-front points, coupling equations of the first-order derivatives and principal second-order derivatives are derived in a dimension-splitting approach by differentiating the jump conditions. Our method obtains second-order accurate solution at the grid points and second-order accurate gradient at the interface. The accurate approximation for the gradient is important in applications where the dynamics of the surface is driven by the jump of the solution gradient at the interface. Our method has more compact finite difference stencils compared with those in CIM2 and is suitable for complex interfaces. We tested our method in three dimensions with complex interfaces, including two protein surfaces, and demonstrated that the solution and the gradient at the interface are uniformly second-order accurate, and the convergence results are very robust. We also tested our method with a moving surface whose normal velocity is given by the jump in the gradient at the interface and achieved second-order accurate interface at the final time.

\section*{Acknowledgement}
This work was funded by NSF Awards 1913144 and 2208465. The authors would like to thank Professor Bo Li for helpful discussions, guidance and support in numerical aspects of the paper. The second author would like to thank Professor Yu-Chen Shu for helpful discussions on CIM.

\bibliographystyle{spphys}
\bibliography{ccim}

\newpage
\section*{Statements and Declarations}

\subsection*{Funding}
This work was funded by NSF Awards 1913144 and 2208465.

\subsection*{Competing Interests}
The authors have no relevant financial or non-financial interests to disclose.

\subsection*{Data Availability}
The code is available at \url{https://github.com/Rayzhangzirui/ccim}.

\newpage
\appendix
\label{appendix}
\section{Differentiation of the jump condition}\label{ss:derivation}
In this appendix, we detail the calculation of the formula for approximation of 
$[\pdv*[2]{u}{x_k}]$, through equations involving terms of $[\nabla^2 u]$, as found in 
Section~\ref{ss:jumpd2u}. 
In the following derivation, quantities related to $f$, $a$, $\epsilon$ and $\tau$ are all known. Our final goal is to write the jump of the second derivatives $[\hess u]$ in terms of the known quantities and the one-sided derivatives $\grad u^-$, $\hess u^-$. Since our interface is smooth, we can consider any smooth extension of $\tau$ and $\sigma$ off the interface, therefore quantities such as $\hess \tau$ and $\grad \tau \cdot \vb n$ are well-defined.  

We first consider jumps of the first derivatives of $u$ at the interface, especially 
in terms of jumps of the normal and tangential derivatives of $u$.  Taking tangential 
derivatives on both sides of $u^+-u^- = [u] = \tau$, we get
$$
[\grad u\cdot \vb s_j] = \grad u^+\cdot \vb s_j-\grad u^-\cdot \vb s_j
 = \grad\tau\cdot \vb s_j,
$$
for $j = 1,\dots,d-1$.  On the other hand, with
\begin{equation*}
\left[\epsilon v\right] 
		=  \epsilon^+[v]+[\epsilon]v^-, \tag{\ref{eq:trick} revisited}
\end{equation*}
we can get, when $v = \nabla u\cdot n$, 
\begin{equation*}
	[\grad u \cdot \vb n] = \frac{1}{\epsilon^+} (\sigma - [\epsilon] (\grad u^- \cdot \vb n)). \tag{\ref{eq:jumpDun} revisited}
\end{equation*}

These equations for the jumps of the normal and tangential derivatives of $u$ can then 
be used to get $[\nabla u]$, from
\begin{equation*}
	[\grad u] = [\grad u \cdot \vb n] \vb n + \sum_{j=1}^{d-1}(\grad \tau \cdot \vb s_j) \vb s_j.   \tag{\ref{eq:jumpgradu} revisited}
\end{equation*}
Thus having handled jumps of first derivatives of $u$, we now turn our attention jumps 
of second derivatives and derive equations for the terms of $[\nabla^2 u]$ in three
ways.

\textbf{Tangential derivative of jump of tangential derivative:}
We can get equations on terms of $[\nabla^2 u]$ by starting with the jump of the 
tangential derivative of $u$ along $\vb s_m$, namely $[\nabla u\cdot s_m]$, and taking 
its tangential derivative along $\vb s_n$, giving $\nabla[\nabla u\cdot s_m]\cdot s_n$.  
This quantity can be written as
\begin{equation*}
	\grad \jump{\grad u \cdot \vb s_m} \cdot \vb s_n	
 = 
	\grad (\grad \tau \cdot \vb s_m) \cdot \vb s_n.
\end{equation*}
However, we also have, using \eqref{eq:jumpgradu}, that
\begin{eqnarray*}
		\grad \jump{\grad u \cdot \vb s_m} \cdot \vb s_n
		& = & \vb s_n^T \jump{\grad^2 u} \vb s_m
		+ \vb s_n^T\grad \vb s_m \jump{\grad u}\\
		& = & \vb s_n^T \jump{\grad^2 u} \vb s_m
		+ \frac{1}{\epsp} \left(\sigma - \jump{\epsilon}\grad u^-\cdot \vb n\right) \vb s_n^T\grad \vb s_m \vb n 
		+ \sum_{j=1}^{d-1} (\grad \tau \cdot \vb s_j) \vb s_n^T \grad \vb s_m \vb s_j
\end{eqnarray*}
and, additionally, that
\begin{eqnarray*}
	\grad (\grad \tau \cdot \vb s_m) \cdot \vb s_n & = & \vb s_n^T \grad^2 \tau \vb s_m + s_n^T \grad \vb s_m \grad \tau\\	
	& = & \vb s_n^T \grad^2 \tau \vb s_m + (\grad \tau \cdot \vb n) \vb s_n^T \grad \vb s_m \vb n + \sum_{j=1}^{d-1} (\grad \tau \cdot \vb s_j) \vb s_n^T \grad \vb s_m \vb s_j.
\end{eqnarray*}
Thus, equating these and solving for the term with jumps in second derivatives of $u$, 
especially simplifying using the fact that
\begin{equation*}
	\grad \vb n~ \vb s_j = - \grad \vb s_j \vb n, \quad j = 1, \dots, d-1,
\end{equation*}
from taking the gradient on both sides of 
$\vb n \cdot \vb s_j = 0$ for $j = 1, \dots, d-1$, we get \eqref{eq:dtau}:
\begin{equation*}
	\vb{s}_n^T [\grad^2u]\vb{s}_m = \vb{s}_n^T \grad^2 \tau \vb{s}_m 
	- \frac{1}{\epsilon^+}( \sigma - [\epsilon] \grad u^- \cdot \vb n) \vb{s}_n^T \grad \vn \vb{s}_m 
	- (\grad \tau \cdot \vn ) \vb{s}_n^T \grad \vn \vb{s}_m.
\end{equation*}

\textbf{Tangential derivative of flux jump:}
For more equations on terms of $[\nabla^2 u]$, we consider the tangential derivative of 
the jump $[\epsilon\grad u\cdot \vb n]$ along $\vb s_m$, which satisfies
$$
\nabla[\epsilon\grad u\cdot \vb n]\cdot \vb s_m = \grad\sigma\cdot\vb s_m,
$$
where $\sigma$ is the jump of the flux.  Expanding, we can get
\begin{eqnarray*}
	\grad \jump{\epsilon \grad u \cdot \vb n} \cdot \vb s_m
	& = & \vb s_m^T \grad (\epsp \jump{\grad u} + \jump{\epsilon} \grad u^-) \vb n + \vb s_m^T \grad \vb n (\epsp \jump{\grad u} + \jump{\epsilon} \grad u^-)\\
	& = & \epsp \vb s_m^T \jump{\hess u} \vb n 
	+ \vb s_m^T \grad \epsp [\grad u] \cdot \vb n 
	+ \jump{\epsilon} \vb s_m^T \hess u^- \vb n + 
	\vb s_m^T \jump{\grad \epsilon} \grad u^- \cdot \vb n \\
	& & + \epsp \vb s_m^T \grad \vb n \jump{\grad u} + \jump{\epsilon} \vb s_m^T \grad \vb n \grad u^-.
\end{eqnarray*}
Now, substituting $[\grad u]$ by equation \eqref{eq:jumpgradu} and rearranging, 
we get \eqref{eq:dsigma}:
\begin{equation*}
	\begin{aligned}
	\vb{s}_m^T [\grad^2u]\vn 
	&= \frac{1}{\epsilon^+} \grad \sigma \cdot \vb{s}_m 
	- \frac{[\epsilon]}{\epsilon^+} \vb{s}_m^T \hess u^- \vn 
	- \frac{[\epsilon]}{\epsilon^+} \vb{s}_m^T \grad \vn \grad u^- \\
	&- \sum_{k=1}^{d-1} (\grad\tau\cdot \vb{s}_k)  \vb{s}_m^T \grad \vn \vb{s}_k 
	- \frac{1}{(\epsilon^+)^2} (\grad \epsilon^+ \cdot \vb{s}_m ) ( \sigma - [\epsilon] \grad u^- \cdot \vb n)\\
	&- \frac{1}{\epsilon^+}  [\grad \epsilon \cdot \vb{s}_m]  (\grad u^- \cdot \vb n).
	\end{aligned}
\end{equation*}

\textbf{Jump of PDE:}
For our final equations on terms of $[\nabla^2 u]$, we consider the original PDE.  
In $\omgp$, we have
\begin{equation*}
	-\epsp \Delta u^+ - \grad \epsp \cdot \grad u^+ + a^+ u^+ = f^+,
\end{equation*}
while in $\omgm$, we have
\begin{equation*}
	-\epsm \Delta u^- - \grad \epsm \cdot \grad u^- + a^- u^- = f^-.
\end{equation*}
Dividing these equations by $\epsp$ and $\epsm$, respectively, and finding the jump from
their difference, we get
\begin{eqnarray*}
		[\lap u] & = & -\jump{\frac{f}{\epsilon}} - \jump{\frac{\grad \epsilon}{\epsilon}\cdot \grad u} + \jump{\frac{a}{\epsilon}u}\\
		& = &-\jump{\frac{f}{\epsilon}}
		-\frac{\grad \epsp}{\epsp} \cdot \jump{\grad u}	
		-\jump{\frac{\grad \epsilon}{\epsilon}}\cdot \grad u^-
		+\frac{a^+}{\epsp}\jump{u}
		+\jump{\frac{a}{\epsilon}}u^-\\
		& = & -\jump{\frac{f}{\epsilon}}
		-\jump{\frac{\grad \epsilon}{\epsilon}}\cdot \grad u^-
		+ \frac{a^+}{\epsp}\jump{u}
		+ \jump{\frac{a}{\epsilon}}u^-\\
		& & -\frac{1}{\epsp}\left( \jump{\grad u \cdot \vb n} (\grad \epsp \cdot \vb n)
		+\sum_{j=1}^{d-1}(\grad \tau \cdot \vb s_j) (\grad \epsp \cdot \vb s_j) \right),
\end{eqnarray*}
which is \eqref{eq:df}.


The left hand side of equations \eqref{eq:dtau} \eqref{eq:dsigma} and \eqref{eq:df}:
\begin{equation*}
	\vb{s}_n^T [\grad^2u]\vb{s}_m,\quad \vb{s}_m^T [\grad^2u]\vn, \quad [\lap u] \quad \text{for $m = 1,\dots,d-1$ and $n = m,\dots,d-1$}
\end{equation*}
where $[\grad^2u]$ are the unknown quantities,
can be written as in the form of a matrix-vector product 
\begin{equation*}
	G \left(\left[\pdv{u}{x_k}{x_l}\right]\right)_{1 \leq k \leq l \leq d}.
\end{equation*}
where $G$ is a matrix that only depends on the normal and the tangent vectors, and the vector is a half-vectorization of the jump of the symmetric Hessian matrix $[\hess u]$.

We can show that the absolute value of the determinant of $G$ is 1 in two and three dimensions.
Since the equations are obtained at some interface point $\hx$, we can use a local coordinate system such that $\vb s_i = \vb e_i$ for $i=1,\dots,d-1$
and $\vb n = \vb e_d$.
By choosing a specific ordering for the equations and the half-vectorization, we can write the matrix-vector product in 2D as 
\begin{equation*}
\begin{pmatrix}
\vb{s}_1 [\grad^2u]\vb{s}_1\\
\vb{s}_1 [\grad^2u] \vb n \\
[\lap u]
\end{pmatrix} = 
\begin{pmatrix}
	1 & 0 & 0\\
	1 & 1 & 0\\
	0 & 0 & 1\\
\end{pmatrix}
\begin{pmatrix}
[u_{xx}]\\
[u_{xy}]\\
[u_{yy}]\\
\end{pmatrix}
\end{equation*}
And in 3D
\begin{equation*}
\begin{pmatrix}
\vb{s}_1 [\grad^2u]\vb{s}_1\\
\vb{s}_2 [\grad^2u]\vb{s}_2\\
[\lap u]\\
\vb{s}_1 [\grad^2u] \vb{s}_2\\
\vb{s}_1 [\grad^2u] \vn\\
\vb{s}_2 [\grad^2u] \vn\\
\end{pmatrix} = 
\begin{pmatrix}
	1 & 0 & 0 & 0 & 0 & 0 \\
	0 & 1 & 0 & 0 & 0 & 0 \\
	1 & 1 & 1 & 0 & 0 & 0 \\
	0 & 0 & 0 & 1 & 0 & 0 \\
	0 & 0 & 0 & 0 & 1 & 0 \\
	0 & 0 & 0 & 0 & 0 & 1 \\
\end{pmatrix}
\begin{pmatrix}
[u_{xx}]\\
[u_{xy}]\\
[u_{yy}]\\
[u_{xy}]\\
[u_{xz}]\\
[u_{yz}]\\
\end{pmatrix}
\end{equation*}
Therefore, the determinant of $G$ is $\pm 1$, depending on the ordering of the equations and the half-vectorization.

\section{High Contrast Problem}\label{ss:highcontrast}
One difficulty of interface problems is the so-called large contrast problem, where the ratio of the coefficients $\epsp/\epsm\gg1$.
There are several works that analyze the high-contrast problems in the context of unfitted Nitsche finite element method \cite{burmanRobustFluxError2018} and unfitted finite element methods \cite{chuNewMultiscaleFinite2010,guzmanFiniteElementMethod2017a,burmanNumericalApproximationLarge2012,burmanRobustFluxError2018}. For example, it can be shown that the flux error estimate is independent of the contrast for a class of unfitted Nitsche finite element methods  \cite{burmanRobustFluxError2018}.

\begin{figure}[!htbp]
\centering
\includegraphics[width=\textwidth]{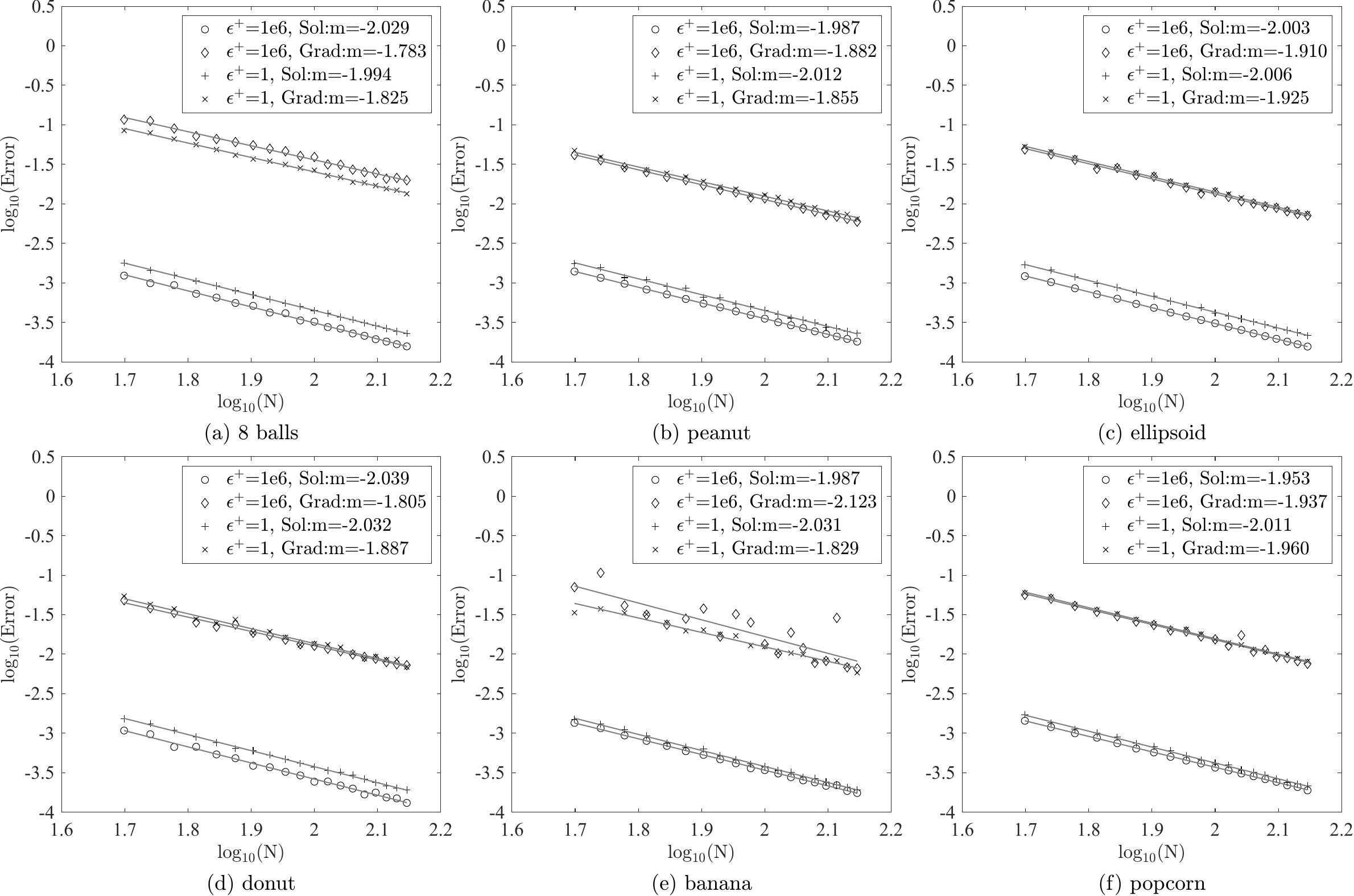}
\caption{The log-log plot of the error versus $N$ for the six surfaces. $\epsm=1$ and $\epsp=1$ or $\epsp=1e6$. In each figure, N ranges from 50 to 140 with the increment $\Delta N = 5$. 
``Sol'' denotes maximum errors of the solution $\norm{u_e-u}_\infty$. 
``Grad'' denotes the maximum errors of the gradient at interface $\norm{\grad u_e - \grad u}_{\infty,\Gamma}$.
$m$ is the slope of the fitting line.}
\label{f:constrast}
\end{figure}

Here we numerically demonstrate that our method is robust for high contrast problems. We consider the same exact solution \eqref{eq:testu} and coefficients \eqref{eq:epsilon} as in Example 1, but with $\epsm=1$ and $\epsp=1$ or $\epsp=1e6$.
Fig \ref{f:constrast} shows the convergence result of the six interfaces.
We see that both the solution and the gradient at the interface are uniformly second-order accurate for both cases.
We also see that the error between the two cases is similar, demonstrating the robustness of our method for high contrast problems.
The theoretical analysis of the high contrast problem is beyond the scope of this paper and will be studied in the future.

\end{document}